\documentclass[3p, 11pt]{elsarticle}

\usepackage{xcolor}
\usepackage{xspace}
\usepackage{graphicx}
\usepackage{setspace}
\usepackage{subcaption}
\usepackage[utf8]{inputenc}
\usepackage{amsmath,amssymb,amsthm,mathrsfs}
\usepackage{lineno}
\usepackage{siunitx}
\usepackage[linesnumbered, ruled, vlined]{algorithm2e}
\usepackage[colorlinks, allcolors=blue, pdfauthor=author,backref=page]{hyperref}
\usepackage[capitalize]{cleveref}
\newcommand{\bm}[1]{\ensuremath{\mathbf{#1}}}
\newcommand{\bs}[1]{\ensuremath{\boldsymbol{#1}}}

\usepackage[english]{babel}
\newtheorem{theorem}{Theorem}

\newtheorem{definition}{Definition}
\newtheorem{remark}{Remark}
\newcommand{\inner}[2]{\big< #1 ,  #2 \big>}
\newcommand{\Inner}[2]{\bigg< #1 ,  #2 \bigg>}

\makeatletter
\def\ps@pprintTitle{%
  \let\@oddhead\@empty
  \let\@evenhead\@empty
  \def\@oddfoot{\reset@font\hfil\thepage\hfil}
  \let\@evenfoot\@oddfoot
}

\makeatother
\begin{document}
\begin{frontmatter}

\title{Adaptive sparse interpolation for accelerating  nonlinear stochastic reduced-order modeling with time-dependent bases}
\author{Mohammad Hossein Naderi}
\author{Hessam Babaee\corref{cor}}
\ead{h.babaee@pitt.edu}
\cortext[cor]{Corresponding author}
\address{Department of Mechanical Engineering and Materials Science, University of Pittsburgh, 3700 O’Hara Street, Pittsburgh, PA 15213, USA,}

\begin{abstract}
    Stochastic reduced-order modeling based on time-dependent bases (TDBs) has proven successful for extracting and exploiting low-dimensional manifold from stochastic partial differential equations (SPDEs). The nominal computational cost of solving a rank-$r$ reduced-order model (ROM) based on time-dependent basis, a.k.a. TDB-ROM, is roughly equal to that of solving the full-order model for $r$ random samples. As of now, this nominal performance can only be achieved for linear or quadratic SPDEs -- at the expense of a highly intrusive process.  On the other hand, for problems with non-polynomial nonlinearity, the computational cost of solving the TDB evolution equations is the same as solving the full-order model.
    In this work, we present an adaptive sparse interpolation algorithm that enables stochastic TDB-ROMs to achieve nominal computational cost for generic nonlinear SPDEs. Our algorithm constructs a low-rank approximation for  the right hand side of the SPDE using the discrete empirical interpolation method (DEIM).
    The presented algorithm does not require any offline computation and as a result the low-rank approximation can adapt to any transient changes of the dynamics on the fly. We also propose a rank-adaptive strategy to control the error of the sparse interpolation. Our algorithm achieves computational speedup by adaptive sampling of the state and random spaces. We illustrate the efficiency of our approach for two test cases: (1) one-dimensional stochastic Burgers’ equation, and (2) two-dimensional compressible Navier-Stokes equations subject to one-hundred-dimensional random perturbations. In all cases, the presented algorithm results in orders of magnitude reduction in the computational cost.
\end{abstract}
\end{frontmatter}

\section{\label{sec:Intro}Introduction}

Propagating uncertainty in evolutionary systems is of major interest to vast applications in science and engineering. However, when the system of interest is governed by a high-dimensional nonlinear dynamical system and when the number of random parameters is large, the computational cost of propagating uncertainty becomes prohibitive. Sampling techniques such as Monte Carlo \cite{kuo_quasi-monte_2012, giles_multilevel_2008} are too expensive and approaches based on polynomial chaos \cite{wiener_homogeneous_1938} suffer from the curse of dimensionality. Extracting and exploiting correlated structures are key ingredients of methods that can significantly reduce the computational cost of solving these problems. Utilizing reduced-order models (ROMs) is one way to exploit the structures, in which the rank of the ROM does not grow exponentially with the dimension of the random space. Instead, the rank of ROM is tied to the intrinsic dimensionality of the system.

The majority of ROM methodologies require an \emph{offline} process for extracting low-rank subspace or manifold from the data. One may have to pay a significant cost in the offline stage with the hope that cost of \emph{online} calculations is much less than solving the full-order model (FOM).  Examples that follow this workflow are  ROMs based on proper orthogonal decomposition (POD)  \cite{sirovich_turbulence_1987, benner_model_2017, schmidt_guide_2020, schmid_dynamic_2010, kutz_dynamic_2016, khodkar_data-driven_2018, Taira_2017, Taira_2020}, in which a low-rank static subspace is extracted by performing a singular value decomposition (SVD) of the matrix of snapshots in the offline stage and an low-order model is built for fast online calculations.  Methods based on deep convolutional autoencoder/decoder also follow the same workflow in which a nonlinear manifold is extracted from   data  \cite{lee_model_2020}. One of the limitations  of this offline-online workflow is the problem of extrapolation to unseen conditions. For example, in POD-ROM, if ROM is utilized in operating conditions (e.g., different Reynolds number, Mach number, boundary condition, etc.) that are different than those conditions that the POD modes are built for, in general, no guarantee can be made about the accuracy of the ROM. This limitation has motivated using \emph{on-the-fly} ROMs, in which the offline stage is eliminated  and the extraction of the low-rank structures as well as building ROM are carried out online. ROMs based time-dependent bases (TDB), a.k.a. TDB-ROMs belong to this category, in which the correlated structures are expressed in the form of a time-dependent subspace.  TDB-ROMs have another advantage to POD-ROMs: Many systems can be very low-dimensional in TDB but very high-dimensional in static bases (POD or dynamic mode decomposition (DMD)). For example, advection-dominated problems with slowly decaying Kolmogorov bandwidth are high-dimensional in POD bases but low-dimensional in TDB \cite{Babaee_2019,kutz_dynamic_2016,Ohlberger_2015}. 

In the context of stochastic reduced-order modeling, dynamically orthogonal (DO) decomposition is the first TDB-ROM methodology for solving stochastic partial differential equations (SPDEs) \cite{Sapsis_Lermusiaux_2009}.
In this method, the random field is decomposed to a set of orthonormal TDBs and their time-dependent stochastic coefficients \cite{Babaee_Choi_Sapsis_Karniadakis_2017}.  Deterministic PDEs are then derived for the evolution of TDBs and the stochastic ROM. Later, it was shown that the DO evolution equations can be obtained from a variational principle: The TDBs and the stochastic coefficients evolve optimally to minimize the residual of the SPDE \cite{Babaee_2019}. It has also been demonstrated the DO closely approximates the instantaneous Karhunen-L\'{o}eve (KL) decomposition of the random field.  Since the introduction of DO, other TDB-ROMs have been proposed. Bi-orthogonal (BO) decomposition \cite{Cheng_Hou_Zhang_2013a} and dynamically bi-orthonormal (DBO) decomposition \cite{patil2020} are two examples of these modifications. These three methodologies (DO, BO, DBO) are all equivalent \cite{choi_equivalence_2014}. They all extract identical low-rank subspaces and they only differ in an in-subspace scaling and rotation. In this study, we utilize the DBO decomposition since the efficiency of this approach for quantifying the uncertainty of highly ill-conditioned physical systems has been established by multiple research studies \cite{patil2020,ramezanianfly_2021, patil_reduced_2021, aitzhan_reduced_2022}. Extracting correlated structures using TDB  has been established in the chemical physics literature for solving high-dimensional deterministic problems before the application of TDB for solving SPDE. In quantum chemistry literature,  the minimization principle whose optimality conditions lead to the evolution equations of TDB is known as the Dirac–Frenkel time-dependent variational principle \cite{beck_multiconfiguration_2000}. Dynamical low-rank approximation \cite{koch_dynamical_2007} also uses the same variational principle for solving matrix differential equations. For the case where the mean flow is not explicitly solved for, the DBO evolution equations in the semi-discrete form (discretized in the physical and random spaces) are the same as the evolution equations obtained for the dynamical low-rank approximation.  

On the other hand, nonlinear POD-based ROMs have the clear advantage that their online evaluation cost is $\mathcal{O}(r)$, where $r$ is the rank of the POD subspace. The low-computational cost of solving POD-ROMs has been made possible due to the recent advances in sparse interpolation and hyper-reduction techniques. One of the most widely used algorithms  is the discrete empirical interpolation method (DEIM) \cite{chaturantabut_nonlinear_2010}, in which the nonlinear terms are discretely sampled at $\mathcal{O}(r)$ points in online calculations. Other methods that aim to extend this idea include the Q-DEIM method \cite{drmac_new_2016}, the Weighted DEIM (W-DEIM) \cite{Drm2018}, Nonlinear DEIM (NLDEIM) \cite{Otto_2019}, and Randomized DEIM (R-DEIM) \cite{saibaba_randomized_2020}.
The localized discrete empirical interpolation method \cite{Peherstorfer_Butnaru_Willcox_Bungartz_2014} was introduced to calculate a number of local subspaces, each tailored to a special part of a dynamical system. 
An adaptivity procedure was introduced in \cite{Peherstorfer_Willcox_2015}, in which the low-rank subspace is updated via an online DEIM sampling strategy. 
The gappy proper orthogonal decomposition (Gappy POD) \cite{everson_karhunenloeve_1995,venturi_gappy_2004} is another approximation technique that was proposed to estimate the nonlinear term with regression instead of interpolation through oversampling. Hyper-reduction techniques approximate the projection of the nonlinear terms onto the POD subspace via sparse interpolation \cite{ryckelynck_priori_2005,Farhat_2015, hernandez_high-performance_2014, antil_two-step_2013}. These approximation techniques for the nonlinear term have been applied successfully to diverse applications. See for example  \cite{hernandez_dimensional_2017, lee_model_2020, chen_eim-degradation_2021, kim_fast_2022, zucatti_calibration_2021, loiseau_noack_brunton_2018, pando_nonlinear_2016}.

Despite the remarkable promise that TDB-ROMs offer for solving high-dimensional problems, the computational cost of solving the TDB-ROM evolution equations precludes the application of these techniques to diverse SPDEs. For certain types of SPDEs, it is possible to achieve significant speedup by using TDB-ROM. However, the speedup is achieved at the cost of a highly intrusive process --- for the derivation and implementation of TDB-ROM equations.   
Assuming that the discretized SPDE has $n$ degrees of freedom in the state space and $s$ random samples are required to achieve statistical convergence, the computational cost of solving full-order model (FOM) scales with $\mathcal{O}(sn^{\alpha})$, where $\alpha \geq 1$ and the value of $\alpha$ depends on the type of SPDE and the spatial discretization. Solving DO, BO and DBO evolution equations for SPDEs with non-polynomial nonlinearities is also  $\mathcal{O}(sn^{\alpha})$, i.e., as expensive as solving FOM. However, the potential, and the promise, of the TDB-ROM is to reduce this cost to $\mathcal{O}(r^2(s+n^{\alpha}))$, where the cost of evolving the TDB and the stochastic coefficients scales with $\mathcal{O}(n^{\alpha})$ and $\mathcal{O}(s)$, respectively and $r$ is the rank of the TDB. Obviously $\mathcal{O}(sn^{\alpha})$ precludes the application of TDB for solving SPDEs with non-polynomial nonlinearity, where both $n$ and $s$ are large numbers. For SPDEs with polynomial nonlinearities, the computational cost of TDB-ROM increases exponentially with the polynomial degree. But even that comes at the cost of a highly intrusive process, which involves replacing the TDB decomposition in the SPDE and carefully deriving the right hand side term-by-term. As we show in this paper, even for non-homogeneous linear SPDE subject to high-dimensional stochastic forcing the computational cost of solving TDB-ROM could be prohibitive.    
The issue of cost is the main reason that the practical applications of DO, BO and DBO have been limited to problems with at most quadratic nonlinearities \cite{Sapsis_Lermusiaux_2009,Cheng_Hou_Zhang_2013a,babaee_minimization_2016, patil2020, cao_stochastic_2018, musharbash_dual_2018}.

In this paper we present a sparse interpolation algorithm to reduce the computational cost of solving DO, BO and DBO equations from $\mathcal{O}(sn^{\alpha})$ to $\mathcal{O}(r^2(s+n^{\alpha}))$. The performance of the presented methodology is agnostic to the type of nonlinearity and therefore it enables achieving the nominal speedup for TDB-ROMs ($\mathcal{O}(r^2(s+n^{\alpha}))$) for generic nonlinear SPDEs. Our algorithm is based on a low-rank approximation of the right hand side of the SPDE via a sparse interpolation algorithm. Our algorithm does not require any offline calculation.  We also  employ adaptive (time-dependent) sampling both in the physical space and the high-dimensional random space. To control the error of this approximation, we propose a rank-adaptive strategy where modes are added and removed to maintain the error below some desired threshold value. 

The remainder of this paper is organized as follows. In \cref{sec:Problem} we present the problem definition and the discretization of the full-order model in the physical and random spaces. In \cref{sec:Method}, we briefly review the DBO decomposition method and introduce the sparse TDB-ROM method. Also, we provide an error bound for our proposed method and present a rank-adaptive algorithm in this section. In \cref{sec:DC}, we show some of the results for the 1D stochastic Burgers' equation and 2D stochastic compressible Navier-Stokes equations. Finally, in \cref{sec:Conclusion}, we present the conclusions.


\section{\label{sec:Problem}Problem Description}

\subsection{Definitions and Notation}
We consider a generic nonlinear SPDE defined by:
\begin{equation}
\begin{aligned}
\frac{\partial v(x, t ; \xi)}{\partial t} &= \mathscr{F}(v(x, t ; \xi)), && x \in D,  \\
v\left(x, t_{0} ; \xi\right) &= v_{0}(x ; \xi), && x \in D,  \\
\mathscr{B}(v(x, t ; \xi)) &= g(x, t), && x \in \partial D, 
\end{aligned}
\label{eq:spde}
\end{equation}
 where $\xi=\{\xi_1, \xi_2, \dots, \xi_d\}$ and $d$ is the number of random parameters, $t \in[0, \infty)$ is time, $x$ denotes the spatial coordinates in the physical domain, $D$ is the physical domain, $\partial D$ is the boundary of the physical domain, $\mathscr{F}$ is a nonlinear differential operator, and $\mathscr{B}$ denotes a linear differential operator that acts on the boundary. The focus of this paper is on uncertainty propagation where the source of uncertainty is uncertain parameters whose joint probability density function is denoted by $\rho(\xi)$. The solution of the above SPDE is denoted by $v(x,t;\xi)$, which is a time-dependent random field. We also denote  the expectation operator with $\mathbb{E}[\sim]$ defined as:
 \begin{equation}
     \mathbb{E}[v(x,t;\xi)] = \int_{\xi} v(x,t;\xi) \rho(\xi) d\xi.
 \end{equation}
 
We denote continuous variables/functions with lower case  ($v$), and we use bold lower case   for vectors ($\bm{v}$), and we use bold upper case  for matrices  ($\bm{V}$).  We use MATLAB convention to indicate elements of a matrix or a vector. For example, if $\bm{p}=[p_1,p_2, \dots, p_r]$ is a set of  indexes, where $p_i$'s are integers in the range of $1\leq p_i \leq n $, then $\bm{V}(\bm{p},:)$ is an $r\times s$ matrix containing the $p_i^{th}$ ($i=1, \dots, r)$  rows of $\bm{V}$.   Similarly, if $\bm{q}= [q_1,q_2, \dots, q_c]$, where $q_i$'s are integers in the range of $1\leq q_i \leq s $, then  $\bm{V}(\bm{p},\bm{q})$ is an $r\times c$ matrix containing the entries of $\bm{V}$ at rows $\bm{p}$ and columns $\bm{q}$.

\subsection{Discretization in the Physical and Random Domains}
The time-dependent field random field $v(x,t;\xi)$ can  be approximated using the following modal decomposition
\begin{equation}
    v(x,t;\xi) \approx \sum_{i=1}^n \sum_{j=1}^s \hat{v}_{ij}(t) \psi_j(\xi) \phi_i(x),
\end{equation}
where, $\phi_i(x)$ are the trial basis functions in the physical domain,  $\psi_j(\xi)$  represent the trial basis functions in the random  space and $\hat{v}_{ij}$ represent the modal coefficients. Let $\bm{x} = [x_1, x_2, \dots, x_n]$ denote the quadrature points in the physical domain and $\bs{\xi}=\{\xi^{(1)}, \xi^{(2)}, \dots, \xi^{(s)} \}$ denote the quadrature points in the random parametric space. Let $\bm{\Phi}=[\phi_1(\bm{x}) | \phi_2(\bm{x}) | \dots | \phi_n(\bm{x})] \in \mathbb{R}^{n\times n}$ be the matrix of basis functions in the physical space evaluated at the quadrature points ($\bm{x}$) and similarly, let 
$\bs{\Psi} = [\psi_1(\bs{\xi}) | \psi_2(\bs{\xi}) | \dots | \psi_s(\bs{\xi})] \in \mathbb{R}^{s \times s}$ denote the matrix of parametric basis functions evaluated at $\bs{\xi}$. Therefore, any deterministic spatial function $u(x)$ evaluated at the quadrature points, i.e., $\bm{u} = u(\bm{x}) \in \mathbb{R}^{n \times 1}$, can be represented via $\bm{u} = \bm{\Phi}\hat{\bm{u}}$ and similarly any function of random parameters $y(\xi)$ evaluated at the quadrature points , i.e., $\bm{y} = y(\bs{\xi}) \in \mathbb{R}^{s \times 1}$, can be expressed via $\bm{y}=\bs{\Psi}\hat{\bm{y}}$. A time-dependent random field  $v(x,t;\xi)$ evaluated at the quadrature points in the physical and random spaces can be  represented via $\bm{V}(t) \in \mathbb{R}^{n\times s}$:  
\begin{equation}
\mathbf{V}(t) =\left[\mathbf{v}_{1}(t) \quad \mathbf{v}_{2}(t)  \,\, \ldots \,\, \mathbf{v}_{s}(t)\right] \in \mathbb{R}^{n \times s},
\end{equation}
where $\bm{v}_i(t)=v(\bm{x},t;\xi^{(i)})$.   
It is straightforward to show that $\bm{V}(t) = \bm{\Phi} \hat{\bm{V}}(t) \bs{\Psi}^T$, where 
 $\hat{\bm{V}}(t) \in \mathbb{R}^{n\times s}$ is the matrix of the modal coefficients, i.e., $\hat{\bm{V}}_{ij}(t)  \equiv \hat{v}_{ij}(t)$.

For the sake of simplicity in our exposition, we consider collocation schemes for discretization in both physical and random spaces, where $\phi_i(x_j) = \delta_{ij}$ and $\psi_i(\xi^{(j)}) = \delta_{ij}$, or alternatively, $\bs{\Phi}=\bm{I}$ and $\bs{\Psi}=\bm{I}$.  As a result, the value of the functions evaluated at the quadrature points is equal to the vector of modal coefficients, i.e., $\bm{u} = \bm{\Phi}\hat{\bm{u}} = \hat{\bm{u}}$ and  $\bm{y} = \bs{\Psi}\hat{\bm{y}} = \hat{\bm{y}}$.  Moreover, $\bm{V}(t)=\hat{\bm{V}}(t)$. 

When a collocation scheme in the random space is used, \cref{eq:spde} can be discretized and solved for each random collocation point independently.  To this end, let the semi-discrete form of \cref{eq:spde} for the collocation point $\xi^{(i)}$ be expressed as: 
\begin{equation}
    \dot{\hat{\mathbf{v}}}_i=\mathcal{F}(\hat{\mathbf{v}}_i), \quad i=1,2,\dots s,
    \label{semi-discrete_pde}
\end{equation}
where $\hat{\bm{v}}_i \in \mathbb{R}^{n}$ is the $( \dot{\sim} ) = d(\sim)/dt $ and $\mathcal{F}( \ . \ )$ represents the discrete representation of  $\mathscr{F}( \ . \ )$ such that
\begin{align*}
    \mathcal{F}: \  & \mathbb{R}^{n} \rightarrow \mathbb{R}^{n}, \\
     &\hat{\bm{v}}_i \rightarrow \mathcal{F}(\hat{\bm{v}}_i).
\end{align*}
Since a collocation scheme is considered for the discretization in the physical domain, in \cref{semi-discrete_pde},  $\hat{\bm{v}}_i$ can be replaced with $\bm{v}_i$. Using the above notation we can write the evolution equations for all samples in the form a \emph{matrix evolution equation} as follows:
\begin{equation}
    \dot{\mathbf{V}}=\mathcal{F}(\mathbf{V}),
    \label{eq:matrix_evol}
\end{equation}
subject to appropriate initial conditions. In \cref{eq:matrix_evol},  when $\mathcal{F}$ is  applied to $\mathbf{V}$, its action is understood to be column wise: $\mathcal{F}(\bm{V})=[\mathcal{F}(\bm{v}_1) \  \mathcal{F}(\bm{v}_2) \  \dots \ \mathcal{F}(\bm{v}_{s}) ] $. Note that the action of $\mathcal{F}$ on each column of $\mathbf{V}$ is independent of the other columns of $\mathbf{V}$. To see how boundary conditions can be incorporated into the above \cref{eq:matrix_evol}, see \cite{Prerna_2021}.

The inner product in both physical and random spaces can be computed using a quadrature rule.  To this end, let $\bm{w}_x = [w_{x_1}, w_{x_2}, \dots, w_{x_n}]$ denote the vector of quadrature weights in the physical domain and  $\bm{w}_{\xi} = [w_{\xi_1}, w_{\xi_2}, \dots, w_{\xi_s}]$ denote the quadrature weights in the random space. Thus, it is possible to approximate the continuous inner product in physical and random spaces as follows: 
\begin{align*}
     \inner{u}{v}_x &= \int_D u(x) v(x) dx \approx \sum_{i=1}^n w_{x_i}u(x_i) v(x_i)=  \bm{u}^T \bm{W}_x \bm{v} ,\\
      \inner{y}{z}_{\xi} &=  \int_{\Omega} y(\xi) z(\xi) \rho(\xi) d\xi \approx \sum_{i=1}^s w_{\xi_i}y(\xi^{(i)})z(\xi^{(i)}) = \bm{y}^T \bm{W}_{\xi} \bm{z},
\end{align*}
where $\bm{W}_x =\mbox{diag}(\bm{w}_x)$ and $\bm{W}_{\xi} =\mbox{diag}(\bm{w}_{\xi})$. Note that in the above definition of the inner product in the random space: $\inner{y}{z}_{\xi}=\mathbb{E}[y(\xi)z(\xi)]$.
Also, a weighted Frobenius norm ($\left\|\mathbf{V}\right\|_{F}$) can be defined as:
\begin{equation}
 \left\|\mathbf{V}\right\|_{F}^{2}=\sum_{i=1}^{n} \sum_{j=1}^{s} \mathbf{W}_{\mathbf{x}_{ii}} \mathbf{W}_{\xi_{jj}} \mathbf{V}_{ij}^{2},
\end{equation}
where $\mathbf{V}_{ij}$ is the entry of matrix $\mathbf{V}$ at the $i^\textit{th}$ row and the $j^\textit{th}$ column. Note that $\left\|\mathbf{V}\right\|_{F}^{2}$ is applied for each instant of time and therefore it is a time-dependent scalar. This norm also approximates:
\begin{equation*}
    \int_D \mathbb{E}[v(x,t;\xi)^2]  dx  \approx \left\|\mathbf{V}\right\|_{F}^{2}.
\end{equation*}

The above setup can be re-purposed for discretizations that  do not have basis functions, for example, finite-difference discretizations in the physical domain or Monte-Carlo sampling in the random space.  For these cases, we can still use $\bm{\Phi}=\bm{I}$ and $\bm{\Psi}=\bm{I}$. For Monte-Carlo sampling, $\bm{W}_{\xi} = \frac{1}{s}\bm{I}$ is used and for finite difference discretizations, $\bm{W}_{x}$ can be taken for example as, $\bm{W}_{x} = \mbox{diag}([\Delta v_1, \Delta v_2, \dots, \Delta v_n])$, where $\Delta v_i$ is the volume of  the cell surrounding the $i^{th}$ grid point. Other higher-order and also non-diagonal weight matrices may be used, for example using the trapezoid rule. 

\subsection{Stochastic Reduced-Order Modeling with TDB}
Our objective is to solve for a low-rank decomposition of $v$ instead of solving \cref{eq:spde}. To this end, we consider the DBO decomposition \cite{patil2020}. As it was shown in \cite{patil2020}, DBO decomposition is equivalent to DO and BO decompositions. However, DBO decomposition shows better numerical performance in comparison to  DO and in contrast to the BO decomposition, the  DBO evolution equations do not become singular when two eigenvalues of the covariance matrix cross. For these reasons, we use DBO decomposition to  demonstrate our methodology. However, as we show in this paper, the  presented algorithm can be  utilized in DO and BO decompositions without any change. The DBO decomposition seeks to approximate $v$ with the following low-rank decomposition: 
\begin{equation}\label{eq:DBO_decomp}
v(x, t; \xi)=\sum_{i=1}^{r} \sum_{j=1}^{r} u_{i}(x, t) \Sigma_{i j}(t) y_{j}(t ; \xi)+e(x, t ; \xi).
\end{equation}
In this representation, $u_{i}(x,t), \, y_{i}(t; \xi), \,\, i=1,2, \ldots, r$ are a set of orthonormal spatial and stochastic modes, respectively:
\begin{equation}
\begin{aligned}
\inner{ u_{i}(x, t)}{ u_{j}(x, t)}_x=\delta_{i j},\\
\inner{y_{i}(t ; \xi)}{ y_{j}(t ; \xi)}_{\xi}=\delta_{i j}.
\end{aligned}
\end{equation}
The stochastic and spatial coefficients are dynamically orthogonal, i.e., the rate of change of these subspaces is orthogonal to the space spanned by these modes:
\begin{equation}
\begin{array}{ll}
\Inner{\dfrac{\partial u_{i}(x, t)}{\partial t}}{ u_{j}(x, t)}_x=0& \quad i, j=1, \ldots, r,\\
\Inner{\dfrac{d y_{i}(t;\xi)}{dt}}{ y_{j}(t ; \xi)}_{\xi}=0& \quad i, j=1, \ldots, r.
\end{array}
\end{equation}

In contrast to the DBO decomposition presented in \cite{patil2020}, in \cref{eq:DBO_decomp}, the mean is not subtracted and as a result $\mathbb{E}[y_i] \neq 0$. This is done for the sake of simplicity, and the presented methodology can be applied to the DBO decomposition where the mean is explicitly subtracted. 
We can write the discrete form of the DBO decomposition as:

\begin{equation}\label{eq:DBO_decomp_dis}
    \mathbf{V}(t) = \mathbf{U}(t) \mathbf{\Sigma}(t)  \mathbf{Y}(t)^{T} + \mathbf{E}(t),
\end{equation}
where,
\begin{equation}
\begin{aligned}
\mathbf{U}(t) &=\left[\mathbf{u}_{1}(t) \quad \mathbf{u}_{2}(t)  \,\, \ldots \,\, \mathbf{u}_{r}(t)\right],\\
\mathbf{Y}(t) &=\left[\mathbf{y}_{1}(t) \quad \mathbf{y}_{2}(t) \,\, \ldots \,\, \mathbf{y}_{r}(t)\right].
\end{aligned}
\end{equation}
Here, $\mathbf{\Sigma}(t)$ is a factorization of the reduced covariance matrix $\mathbf{C}(t) \in \mathbb{R}^{r \times r}$ as in $\mathbf{C}(t)=\mathbf{\Sigma}(t) \mathbf{\Sigma}(t)^T$, and $\mathbf{E}(t)$ is the reduction error. The columns of $\mathbf{U}(t)$ and $\mathbf{Y}(t)$ are a set of orthonormal spatial and stochastic modes, respectively. Given the above definitions the orthonormality of the modes in the discrete form implies that:
\begin{subequations}\label{eq:const}
\begin{align}
\mathbf{U}(t)^T \bm{W}_x \mathbf{U}(t)&= \mathbf{I}, \label{eq:const_U} \\
\mathbf{Y}(t)^T \bm{W}_{\xi} \mathbf{Y}(t)&=\mathbf{I}. \label{eq:const_Y}
\end{align}
\end{subequations}
Here, we introduce the notion of presenting the matrices and equations  in the \emph{compressed} and \emph{decompressed} forms in the context of the TDB decomposition.
\begin{definition}
 Let $\tilde{\bm{V}} = \bm{U}\bm{\Sigma}\bm{Y}^T$ be the DBO decomposition as presented in \cref{eq:DBO_decomp_dis}, then $\tilde{\bm{V}}$ is in the compressed form if it is expressed versus the triplet $\{\bm{U},\bm{\Sigma},\bm{Y}\}$. However, $\tilde{\bm{V}}$ is in the decompressed form if  the factorized matrices are multiplied and $\tilde{\bm{V}}$ is formed explicitly. 
\end{definition}
The memory requirement for storing $\tilde{\bm{V}}$ in the compressed form is $r(n+s+r)$ and in the decompressed form is $sn$.  Obviously for computational purposes it is highly advantageous to avoid decompressing any quantity. The notion of compressed and decompressed can be extended to BO and DO decompositions, which factorize $\tilde{\bm{V}}$ to two matrices: $\tilde{\bm{V}} = \bm{U}_{BO/DO}\bm{Y}^T_{BO/DO}$.   

\subsection{Variational Principle}
The central idea behind reduced-order modeling based on TDBs is that the bases evolve optimally to minimize the ROM residual. The residual is obtained by replacing the DBO decomposition into FOM, given by \cref{semi-discrete_pde}. Because DBO is a low-rank approximation, $\bm{E}(t) \neq \bm{0}$ and therefore the DBO decomposition cannot satisfy the FOM exactly and there will be a residual equal to:
\begin{equation}
    \bm{R} = \frac{d\left(\mathbf{U} \mathbf{\Sigma} \mathbf{Y}^{T}\right)}{d t}-\mathcal{F}(\mathbf{U} \mathbf{\Sigma} \mathbf{Y}^{T}).
\end{equation}
The evolution equations of the DBO components are obtained by minimizing the residual as shown below:
\begin{equation}
\mathcal{G}(\dot{\mathbf{U}}, \dot{\mathbf{\Sigma}}, \dot{\mathbf{Y}})=\left\|\frac{d\left(\mathbf{U} \mathbf{\Sigma} \mathbf{Y}^{T}\right)}{d t}-\mathcal{F}(\mathbf{U} \mathbf{\Sigma} \mathbf{Y}^{T})\right\|_{F}^2,
\label{eq:VP}
\end{equation}
subject to  the orthonormality constraints given by \cref{eq:const}.  
The variational principle aims to minimize the residual by optimally updating $\mathbf{U}$, $\mathbf{\Sigma}$, and $\mathbf{Y}$.  The optimality conditions of the variational principle result in closed-form evolution equations for $\mathbf{U}$, $\mathbf{\Sigma}$, and $\mathbf{Y}$. As indicated in \cite{Prerna_2021}, the closed-form evolution equations of the DBO decomposition are defined by: 
\begin{subequations}\label{eq:DBO_evol}
\begin{align}
&\dot{\mathbf{\Sigma}}=\mathbf{U}^{T} \mathbf{W}_{x} \mathbf{F} \mathbf{W}_{\xi} \mathbf{Y} \label{eq:DBO_evol_S}\\ 
&\dot{\mathbf{U}}=\left(\mathbf{I}-\mathbf{U} \mathbf{U}^{T} \mathbf{W}_{x}\right) \mathbf{F} \mathbf{W}_{\xi} \mathbf{Y} \mathbf{\Sigma}^{-1}, \label{eq:DBO_evol_U} \\  
&\dot{\mathbf{Y}}=\left(\mathbf{I}-\mathbf{Y} \mathbf{Y}^{T} \mathbf{W}_{\xi}\right) \mathbf{F}^{T} \mathbf{W}_{x} \mathbf{U} \mathbf{\Sigma}^{-T}.  \label{eq:DBO_evol_Y}
\end{align}
\end{subequations}
In the above equations, $\mathbf{F} \in \mathbb{R}^{n\times s}$ is defined as $\mathbf{F} =\mathcal{F}(\mathbf{U} \mathbf{\Sigma} \mathbf{Y}^{T})$. The above variational principle is the same as the Dirac–Frenkel time-dependent variational principle in the quantum chemistry literature \cite{beck_multiconfiguration_2000} or the  dynamical low-rank approximation \cite{koch_dynamical_2007}. As it was shown in \cite{Babaee_2019}, it is possible to derive a similar variational principle for the DO decomposition, whose optimality conditions lead to the DO evolution equations.  

\subsection{Computational Cost}
The computational cost of solving \cref{eq:DBO_evol_S,eq:DBO_evol_U,eq:DBO_evol_Y} for a  general nonlinear SPDE scales the same as of that  solving the FOM. The computational cost of computing the right-hand side of the FOM  is $\mathcal{O}(sn^{\alpha})$ for any sampling-based strategy (e.g., Monte-Carlo or probabilistic collocation methods).  Here the computational cost of evaluating $\mathcal{F}(\bm{v})$ is $\mathcal{O}(n^{\alpha})$, where $\alpha \geq 1$ depends of the type of SPDE and the discretization. For sparse spatial discretization schemes,  for example finite difference or spectral element method, and explicit time-advancement  $\alpha=1$. Since we consider only sparse spatial discretizations and explicit time-advancement in this paper, we choose $\alpha=1$ in the rest of this paper. In \cref{eq:DBO_evol_S,eq:DBO_evol_U,eq:DBO_evol_Y}, first $\bm{F}$ needs to be formed at the cost of  $\mathcal{O}(sn)$, which is the same cost of the FOM. In addition to this cost,  in \cref{eq:DBO_evol_U,eq:DBO_evol_Y},  $\bm{F}$ is projected onto $\bm{U}$ and $\bm{Y}$, where each projection requires $\mathcal{O}(rns)$ operations. There are other auxiliary operations, but they are negligible to computing and projecting $\bm{F}$ onto the spatial and stochastic subspaces. 

The DBO equations presented in  \cref{eq:DBO_evol_S,eq:DBO_evol_U,eq:DBO_evol_Y} are in the \emph{decompressed form}, where $\bm{F}$ is formed explicitly.   The source of this difficulty is the nonlinear terms on the right hand side of the SPDE, which requires decompressing $\bm{V}=\bm{U}\bm{\Sigma}\bm{Y}^T$ and applying the nonlinear map $\mathcal{F}$ on every column of  $\bm{V}$. Because $n$ and $s$ are often very large, the computational cost and memory requirement of  $\mathcal{O}(sn)$  is prohibitive for most applications.   This cost can be avoided for simple $\mathcal{F}$, i.e.,   linear  SPDEs and quadratic nonlinearities. For linear and quadratic $\mathcal{F}$, one can plug in the DBO decomposition into $\mathcal{F}$ and derive a compressed form of $\mathcal{F}$ in terms of the DBO components. However,  this process is highly intrusive as it requires derivation and implementation of new evolution equations, which is discussed further in \ref{AP_Computational Cost}. For cubic and higher polynomial nonlinearity,  the number of terms generated via this process grows exponentially fast with the polynomial order. As it is shown in \ref{AP_Computational Cost}, even for non-homogeneous linear SPDEs, the computational cost of computing the right hand side of \cref{eq:DBO_evol_S,eq:DBO_evol_U,eq:DBO_evol_Y} is $\mathcal{O}(rns)$, which is prohibitive for large $n$ and $s$. The DO and BO formulations also have the same computational cost scaling. The presented methodology seeks to approximate $\bm{F}$ in a cost-effective way and although we present our methodology for the DBO formulation,  the algorithm is directly applicable to DO and BO or more broadly TDB-ROM for SPDEs.


\section{\label{sec:Method}Methodology}

We present a  methodology to reduce the computational cost of evaluating the right hand side term to $\mathcal{O}(r^2(n+s))$ for both non-homogeneous linear equations as well as  nonlinear SPDEs with generic nonlinearity (polynomial and non-polynomial).

\subsection{On-the-fly Sparse Interpolation}\label{sec:sparse_interp}
Our approach to reduce the computational cost of solving the TDB-ROM equations is to represent $\bm{F}(t)$ in the compressed form by constructing a low-rank approximation of $\bm{F}(t)$ on the fly. To this end,  we present a sparse interpolation algorithm to interpolate columns of $\bm{F}(t)$ onto a low-rank time-dependent bases denoted by $\bm{U}_F(t)$. We present our algorithm for explicit time integration schemes. The sparse interpolation algorithm for DBO is different from the DEIM algorithm used in POD-ROM in several ways.  First, to preserve the advantages that TDB-ROM offers, the basis for the right hand side of SPDE ($\bm{U}_F$) should be   time-dependent. Moreover, $\bm{U}_F$ ideally needs to be computed on the fly and not from data. Computing $\bm{U}_F$ from data would detract from some of the key advantages of reduced order modeling based on the time-dependent bases. 
Second, in the TDB-ROM formulation, one needs a matrix interpolation algorithm rather than a vector interpolation scheme as is the case in POD-ROM. Third, in the DBO formulation, the computational cost issues are not only limited to the nonlinear terms. Even for linear non-homogeneous SPDEs, the computational cost is still prohibitive.

\begin{figure}[!t]
    \centering
    \includegraphics[width=\textwidth]{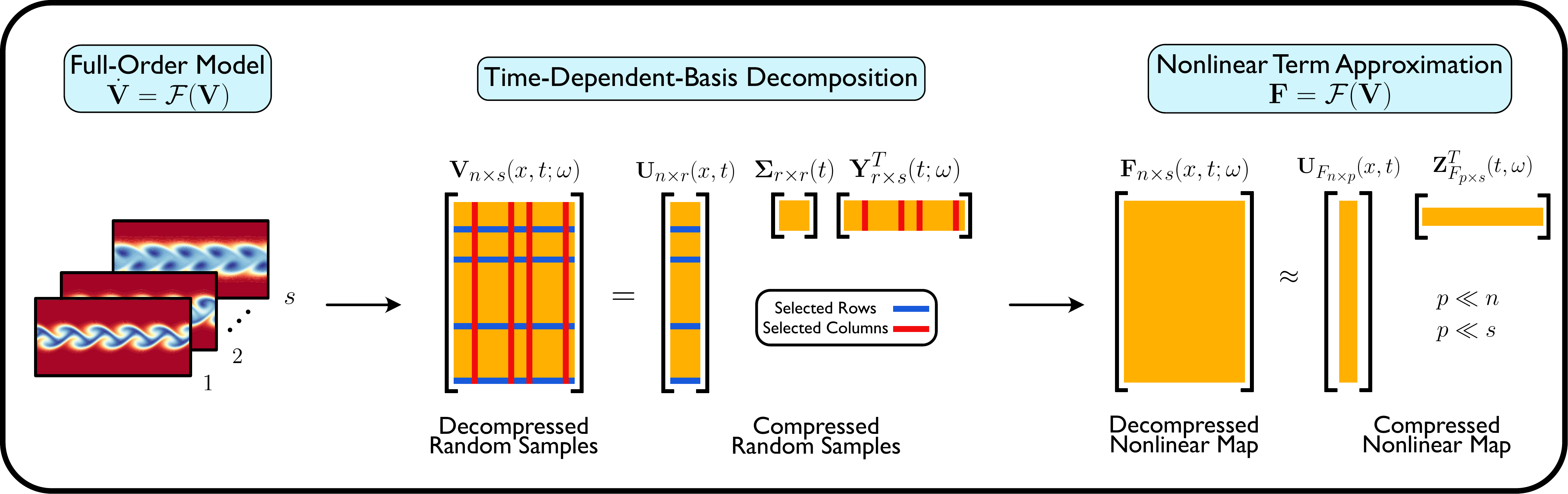}
    \caption{Schematic of the sparse interpolation  algorithm.}
    \label{fig:SDBO}
\end{figure}

The presented algorithm is informed by the above considerations. In particular, we seek to approximate the entire right hand side of the SPDE, i.e., linear and nonlinear terms, with  $\mathbf{F}(t) = \hat{\bm{F}}(t)  + \bm{E}_F(t)$, where:
\begin{equation}
 \hat{\bm{F}}(t) = \bm{U}_F(t) \bm{Z}_F^T(t),
    \label{eq:F}
\end{equation}
and $\bm{U}_F(t) \in \mathbb{R}^{n\times p}$. The columns of $\bm{U}_F(t)$  are  a set of instantaneously orthonormal spatial modes,  i.e., $\bm{U}_F(t)^T \bm{W}_x \bm{U}_F(t) = \bm{I}$,  $\bm{Z}_F(t)  \in \mathbb{R}^{s\times p}$ is the matrix of random coefficients and $\bm{E}_F(t)  \in \mathbb{R}^{n\times s}$ is the approximation error. It is instructive to note that \cref{eq:F}  is the discrete representation of the  continuous form as given below:
\begin{equation}\label{eq:F_cont}
     \mathscr{F}(v(x, t ; \xi)) = \hat{\mathscr{F}}(v(x, t ; \xi))  + e_{\mathscr{F}}(x, t ; \xi), \quad \mbox{where} \quad \hat{\mathscr{F}}(v(x, t ; \xi)) = \sum_{i=1}^p u_{\mathscr{F}_i}(x,t) z_{\mathscr{F}_i}(t;\xi),
\end{equation}
In \cref{eq:F_cont}, $u_{\mathscr{F}_i}(x,t)$,    $z_{\mathscr{F}_i}(t;\xi)$ and $e_{\mathscr{F}}(x, t ; \xi)$ are the continuous representations of $\bm{U}_F(t)$,  $\bm{Z}_F(t)$ and $\bm{E}_F(t)$, respectively.   \cref{eq:F} is an instantaneous low-rank  approximation of matrix $\bm{F}(t)$, where $p \ll n$ and $p \ll s$.   Minimal approximation error is achieved if $\hat{\bm{F}}(t)$ is the  rank-$p$ SVD truncated approximation of $\bm{F}(t)$. We  show that the presented algorithm closely approximates the optimal SVD low-rank approximation. 
In the following, we show how we can compute $\bm{U}_F(t)$ and $\bm{Z}_F(t)$ by sampling $p$ columns and $p$ rows of $\bm{F}(t)$ as shown in \cref{fig:SDBO}.

We present our algorithm for explicit time integration of  \cref{eq:DBO_evol_S}-\cref{eq:DBO_evol_Y} and we  drop the explicit dependence on time for brevity. Instead, we denote quantities from the previous time step  with the superscript $( \sim )^{*}$ and quantities the current time step   are shown without any superscript.  
\\

\noindent \textbf{Computing $\bm{U}_F$:} Ideally $\bm{U}_F$ should be the $p$ dominant left singular vectors of $\bm{F}$. However, computing the exact left singular vectors of $\bm{F}$ requires computing  all entries of matrix $\bm{F}$, which we want to avoid  in the first place. Instead, we seek to approximate the $p$ dominant left singular vectors of $\bm{F}$ by only computing $p$ columns of $\bm{F}$. Then we perform SVD of the obtained matrix.  Let the indexes of the $p$ selected columns be denoted by: $\bm{q} = [q_1,q_2, \dots, q_p]$. Because the columns of $\bm{F}$ can be computed independently from each other,  we can compute $\bm{F}(:,\bm{q})$ by constructing the $p$ samples from the TDB expression: $\bm{F}(:,\bm{q}) = \mathcal{F}(\bm{U}\bm{\Sigma}\bm{Y}(\bm{q},:)^T)$. To compute the SVD we first form the correlation matrix: $\bm{C}_{F} = \bm{F}(:,\bm{q})^T \bm{W}_x \bm{F}(:,\bm{q})$,  where $\bm{C} \in \mathbb{R}^{p\times p}$. Then, we compute the eigenvectors and eigenvalues of the correlations matrix: $\bm{C}_F\bm{\Psi}_F = \bm{\Psi}_F\bm{\Lambda}_F$, where  $\bm{\Psi}_F$ and $\bm{\Lambda}_F$ are the matrices of eigenvectors and eigenvalues of $\bm{C}_F$, respectively. The left singular vectors of $\bm{F}(:,\bm{q})$ are then obtained by:
\begin{equation}\label{eq:eig_F_q}
\bm{U}_F = \bm{F}(:,\bm{q})\bm{\Psi}_F\bm{\Lambda}_F^{-1/2}
\end{equation}
 The computational complexity of the above the computation is $\mathcal{O}(p^2 n)$ -- ignoring the  cost of computing the eigen-decomposition of  $\bm{C}_F$ because $p\ll n$.

Obviously, the choice of selected columns  ($\bm{q}$) is critical to obtain near optimal approximation for $\bm{U}_F$. In this work, we use sparse sampling indexes obtained either from DEIM or Q-DEIM algorithms because of the near optimal performance these two algorithms exhibit.  Both of these algorithms require the left singular vectors of matrix $\bm{F}^T$ or equivalently, the right singular vectors of matrix $\bm{F}$. However, again this is a circular problem: if we somehow had the SVD of $\bm{F}$, we would use the SVD of $\bm{F}$ instead of \cref{eq:F}. Fortunately, we only need the SVD of $\bm{F}$ for selecting which columns to sample and in practice a close approximation of the right singular vectors of $\bm{F}$ is sufficient for column selection. To this end,  we use the right singular vectors of $\hat{\bm{F}}$ from the previous time step. As we show below, we can effectively compute the right singular vectors of $\hat{\bm{F}}^*$   using the low-rank approximation given by \cref{eq:F}. It is straightforward to show that the right singular vectors of $\hat{\bm{F}}^*$ are the same as the left singular vectors of $\bm{Z}_F^*$.   This can be shown by noting that the right singular vectors of $\hat{\bm{F}}^*$ are the eigenvectors of $\hat{\bm{F}}^{*^T} \bm{W}_x  \hat{\bm{F}}^*$ and  the left singular vectors of $\bm{Z}_F^*$ are the eigenvectors of  $\bm{Z}_F^*{\bm{Z}_F^*}^T$. However, $\hat{\bm{F}}^{*^T} \bm{W}_x  \hat{\bm{F}}^*$ and $\bm{Z}_F^*{\bm{Z}_F^*}^T$ are equal to each other:
\begin{equation*}
    \hat{\bm{F}}^{*^T} \bm{W}_x  \hat{\bm{F}}^* = \bm{Z}_F^*{\bm{U}_F^*}^T  \bm{W}_x  \bm{U}^*_F {\bm{Z}_F^*}^T = \bm{Z}_F^* {\bm{Z}_F^*}^T,
\end{equation*}
and therefore, their eigenvectors are identical. Here, we have used the orthonormality of the columns of $\bm{U}^*_F$. Because $\bm{Z}^*_F$ is a tall and skinny matrix ($p\ll s$), to compute the right singular vectors of $\bm{Z}^*_F$, it is better to avoid forming $\bm{Z}^*_F{\bm{Z}^*_F}^T$, because $\bm{Z}^*_F{\bm{Z}_F^*}^T$ is a large matrix and moreover it is rank deficient.  Instead, we  compute $\bm{C}^*_Z = {\bm{Z}_F^*}^T\bm{W}_{\xi}\bm{Z}^*_F \in \mathbb{R}^{p\times p}$, which is a small matrix and full rank and follow these steps:
\begin{equation}\label{eq:Y_f}
\bm{C}^*_Z \bm{\Psi}^*_Z =  \bm{\Psi}^*_Z {\bm{\Sigma}^*_Z}^2 \quad \longrightarrow \quad \bm{Y}^*_F = \bm{Z}^*_F\bm{\Psi}^*_Z{\bm{\Sigma}^*_Z}^{-1},
\end{equation}
where $\bm{Y}^*_F \in \mathbb{R}^{s \times p}$ is the matrix of right singular vectors of $\hat{\bm{F}}^*$ and ${\bm{Y}_F^*}^T \bm{W}_{\xi} \bm{Y}^*_F= \bm{I}$. The matrix $\bm{Y}^*_F$ represents a low-rank subspace  for the rows of $\hat{\bm{F}}^*$ and it can be used in DEIM or Q-DEIM algorithms to obtain the indices of the selected columns ($\bm{q}$). Once the indices of the selected columns are determined, $\bm{U}_F$ (for the current time step) can be computed using  $\bm{U}_F=\bm{F}(:,\bm{q})\bm{\Psi}_F\bm{\Lambda}^{-1/2}$ as explained above. The computational cost of computing $\bm{Y}^*_F$ is $\mathcal{O}(p^2s)$. In order to find $\bm{q}$ in the first time step, we need to know the value of $\mathbf{Z}^*_F$. Note that in the presented algorithm, only the sampling columns are obtained based on the previous time step solution and the values of the columns are computed for the time step in question.

For the first time step,  since  $\mathbf{Z}^*_F$ from the initial condition (at $t=0$) is not available, we use $\mathbf{Y}^*$ as the input of the sparse selection algorithm at step two in \cref{alg:S-TDB-ROM}, which approximately represents the nonlinear basis for finding $\bm{q}$. Next, we follow steps 3-9 in \cref{alg:S-TDB-ROM} to compute the $\mathbf{Z}_F$ and use this matrix at $t=0$.  Also, for explicit time integration schemes that require right-hand side evaluations at substeps, for example, explicit Runge-Kutta schemes, or multi-step time integrators, the indices $\bm{q}$ may be updated at each substep/step or computed once at the first substep/step. We have chosen the latter approach in all examples presented in this paper.   Note that any additional error resulted from the suboptimal determination of $\bm{q}$ ---  either by how $\bm{q}$ is computed  at the initial condition or whether $\bm{q}$ is not updated in each substep/step ---  is controlled by the adaptive algorithm presented in  \cref{sec:adaptive}, where the number of sampled columns increases to meet the reconstruction error threshold.  One can also devise an iterative algorithm to find the optimal $\bm{q}$ ---  perfecting  $\bm{q}$  obtained   from the   previous time step by iteratively repeating steps 1-9 in \cref{alg:S-TDB-ROM}.  However, the cost of an iterative algorithm can easily exceed that of sampling additional suboptimal columns.   \\   

\noindent \textbf{Computing $\bm{Z}_F$:} The columns of $\bm{U}_F$ constitute a low-rank basis for the columns of $\bm{F}$ and they closely approximate the $p$ dominant left singular vectors of $\bm{F}$. We utilize $\bm{U}_F$  and apply DEIM interpolatory projection to all columns of $\bm{F}$. In particular,  we apply DEIM or Q-DEIM algorithms to $\bm{U}_F$ to select $p$ rows. Let $\bm{p}=[p_1,\dots, p_p]$ be the integer vector containing the indices of the selected rows and $\bm{P}=[\bm{e}_{p_1}, \dots, \bm{e}_{p_p}] \in \mathbb{R}^{n \times p}$ is a matrix obtained by selecting certain columns of the identity matrix, where   $\mathbf{e}_{p_i}$ is the $p_i^\textit{th}$ column of the identity matrix.  For example, if $p=3$ and $n=100$ and $\bm{p}=[20, 17, 84]$, then matrix $\bm{P}$ is of size $100 \times 3$ and elements of matrix $\bm{P}$ are all zero except: $\bm{P}(20,1) = \bm{P}(17,2) =\bm{P}(84,3)=1$. Therefore, $\bm{P}^T\bm{F} = \bm{F}(\bm{p},:)$ and  $\bm{P}^T\bm{U}_F \equiv \bm{U}_F(\bm{p},:)$.  Then we sample $\bm{F}$ at the $p$ selected rows, i.e., we compute $\bm{F}(\bm{p},:)$. 

To compute  $\bm{F}(\bm{p},:)$, we need to calculate various spatial derivatives and therefore, we need to know the values of adjacent points. This step depends on the numerical scheme used for the spatial discretization of the SPDE. For example, if we use the spectral element method, to calculate the derivative at a selected spatial point we need to have the values of other points in that element. We denote the index of adjacent points with $\bm{p}_a$. Note that $\bm{F}$ is calculated at points indexed by $\bm{p}$, however the values of $\bm{V}$ at the adjacent points must be provided, which can be obtained via the TDB expansion, i.e., $\hat{\bm{V}}([\bm{p},\bm{p}_a],:) = \bm{U}([\bm{p},\bm{p}_a],:)\bm{\Sigma}\bm{Y}^T$.
After computing $\bm{F}(\bm{p},:)$, the coefficient $\bm{Z}^T_F$ is obtained by:
\begin{equation*}
    \bm{Z}^T_F = \left(\mathbf{P}^{T} \bm{U}_F\right)^{-1} \mathbf{P}^{T} \bm{F} =   \bm{U}_F(\bm{p},:)^{-1}\bm{F}(\bm{p},:).
\end{equation*}
Therefore,
\begin{equation}\label{eq:Z_f}
    \bm{Z}_F = \bm{F}(\bm{p},:)^T   \bm{U}_F(\bm{p},:)^{-T}.
\end{equation}
The computational cost of computing $\bm{Z}_F$ is therefore  $\mathcal{O}(ps)$ due to the computation of right hand side at $p$ points in the physical domain for all $s$ samples.

\begin{algorithm}[t]
\setstretch{1.35}
\SetAlgoLined
\KwIn{$\mathbf{U} \in \mathbb{R}^{n \times r}$, $\mathbf{\Sigma} \in \mathbb{R}^{r \times r}$, $\mathbf{Y}\in \mathbb{R}^{s \times r}$, $\mathbf{Z}^{*}_F \in \mathbb{R}^{s \times p}$, $\mathbf{W}_{x} \in \mathbb{R}^{n \times n}$, $\mathbf{W}_{\xi} \in \mathbb{R}^{s \times s}$ (Data from the previous time step is indicated by superscript $( \sim )^{*}$)}
\KwOut{$\dot{\mathbf{U}}$, $\dot{\mathbf{\Sigma}}$, $\dot{\mathbf{Y}}$, $\mathbf{Z}_F$}
$\mathbf{Y}^{*}_F \leftarrow$ \texttt{SVD(${\mathbf{Z}^{*}_F}^T$)} \hspace{2.535cm} $\rhd$ find left-singular vectors of ${\mathbf{Z}^{*}_F}^T$ ($\mathbf{Y}^{*}_F \in \mathbb{R}^{s \times p}$)\;
${\bm{q}} \leftarrow$ \texttt{Sparse\_Selection($\mathbf{Y}^{*}_F$)} \hspace{0.5cm}  $\rhd$ apply sparse selection algorithm to find $p$ columns\;
$\mathbf{V}_{\bm{q}} = \mathbf{U} \mathbf{\Sigma} \mathbf{Y}(\bm{q},:)^{T}$ \hspace{2.43cm} $\rhd$ construct $\mathbf{V}_{\bm{q}} \in \mathbb{R}^{n \times p}$ with selected columns\;
$\mathbf{F}_{\bm{q}} = \mathcal{F}(\mathbf{V}_{\bm{q}})$ \hspace{3.385cm} $\rhd$ compute $\mathbf{F}_{\bm{q}} \in \mathbb{R}^{n \times p}$\;
$\bm{U}_F \leftarrow$ \texttt{SVD($\mathbf{F}_{\bm{q}}$)} \hspace{2.825cm} $\rhd$ find left-singular vectors of $\mathbf{F}_{\bm{q}}$ ($\bm{U}_F \in \mathbb{R}^{n \times p}$)\;
$\bm{p} \leftarrow$ \texttt{Sparse\_Selection($\bm{U}_F$)} \hspace{0.491cm} $\rhd$ apply sparse selection algorithm to find $p$ rows\;
$\bm{V}([\bm{p},\bm{p}_a],:) = \bm{U}([\bm{p},\bm{p}_a],:)\bm{\Sigma}\bm{Y}^T$ \hspace{0.180cm} $\rhd$ construct $\bm{V}([\bm{p},\bm{p}_a],:) \in \mathbb{R}^{(p+p_a) \times s}$ with selected rows\;
$\bm{F}_{\bm{p}} = \mathcal{F}(\bm{V}([\bm{p},\bm{p}_a],:))$ \hspace{1.96cm} $\rhd$ compute $\mathbf{F}_{\bm{p}} \in \mathbb{R}^{p \times s}$\;
$\bm{Z}_F = \mathbf{F}_{\bm{p}}^T  \big( \bm{U}_F(\bm{p},:)\big)^{-T}$ \hspace{1.775cm} $\rhd$ form $\bm{Z}_F$ with selected rows ($\bm{p}$) of $\bm{U}_{F}$\;
$\dot{\mathbf{\Sigma}}=\big(\mathbf{U}^{T} \mathbf{W}_{x} \bm{U}_F\big) \big(\mathbf{Z}_F^T \mathbf{W}_{\xi} \mathbf{Y}\big)$\;
$\dot{\mathbf{U}}=\big((\mathbf{I}-\mathbf{U} \mathbf{U}^{T} \mathbf{W}_{x}) \bm{U}_F\big) \big(\mathbf{Z}_F^T \mathbf{W}_{\xi} \mathbf{Y} \mathbf{\Sigma}^{-1}\big)$\;
$\dot{\mathbf{Y}}=\big((\mathbf{I}-\mathbf{Y} \mathbf{Y}^{T} \mathbf{W}_{\xi}) \mathbf{Z}_F\big) \big(\bm{U}_F^{T} \mathbf{W}_{x} \mathbf{U} \mathbf{\Sigma}^{-T}\big)$;
\caption{Sparse TDB-ROM Algorithm}
\label{alg:S-TDB-ROM}
\end{algorithm}

Finally, we can replace $\mathbf{F}$ in the closed-form evolution equations of the DBO decomposition by $\bm{U}_F \bm{Z}_F^T$ to obtain the new evolution  equations as in the following:
\begin{subequations}
\begin{align}
\dot{\mathbf{\Sigma}}&=\big(\mathbf{U}^{T} \mathbf{W}_{x} \bm{U}_F\big) \big(\mathbf{Z}_F^T \mathbf{W}_{\xi} \mathbf{Y}\big), \label{eq:SDBO_evol_S}\\
\dot{\mathbf{U}}&=\big((\mathbf{I}-\mathbf{U} \mathbf{U}^{T} \mathbf{W}_{x}) \bm{U}_F\big) \big(\mathbf{Z}_F^T \mathbf{W}_{\xi} \mathbf{Y} \mathbf{\Sigma}^{-1}\big), \label{eq:SDBO_evol_U} \\
\dot{\mathbf{Y}}&=\big((\mathbf{I}-\mathbf{Y} \mathbf{Y}^{T} \mathbf{W}_{\xi}) \mathbf{Z}_F\big) \big(\bm{U}_F^{T} \mathbf{W}_{x} \mathbf{U} \mathbf{\Sigma}^{-T}\big). \label{eq:SDBO_evol_Y}
\end{align}
\end{subequations}

\cref{eq:SDBO_evol_S,eq:SDBO_evol_U,eq:SDBO_evol_Y} are the DBO evolution equation in the compressed form. Solving DBO equations in this form offers significant computational advantages in comparison to \cref{eq:DBO_evol_S,eq:DBO_evol_U,eq:DBO_evol_Y} as enumerated below:
\begin{itemize}
    \item In \cref{eq:SDBO_evol_S,eq:SDBO_evol_U,eq:SDBO_evol_Y}, the matrix $\bm{F}$ is never fully formed and instead $\bm{U}_F$ and $\bm{Z}_F$ are computed.  This reduces the computational cost of $\mathcal{O}(sn)$ to $\mathcal{O}(p^2(s+n))$.  
    \item Solving \cref{eq:DBO_evol_S,eq:DBO_evol_U,eq:DBO_evol_Y} requires either storing the full ($n\times s$) matrix $\bm{F}$ or computing the columns of $\bm{F}$ one (or several) at a time.  On the other hand, solving the TDB-ROM in the form of \cref{eq:SDBO_evol_S,eq:SDBO_evol_U,eq:SDBO_evol_Y} does not require storing any matrix of size $n\times s$ and the memory requirement is reduced to $p(s+n)$. 
    \item  Solving \cref{eq:SDBO_evol_S,eq:SDBO_evol_U,eq:SDBO_evol_Y} is significantly less intrusive than solving the DBO equations in the compressed form. Computing $\bm{U}_F$ requires computing the right hand side of the FOM for $p$ samples, which can be done in a black-box fashion for many solvers. Computing $\bm{Z}_F$ requires knowing the governing equation and sampling the right hand side at $p$ grid points.  However, solving \cref{eq:SDBO_evol_S,eq:SDBO_evol_U,eq:SDBO_evol_Y} does not require replacing the DBO decomposition on the right hand side and working out the expansion term by term, as it is  done for linear and quadratic nonlinear SPDEs. 
\end{itemize}
We also note that one of the by-products of the presented methodology is an efficient algorithm for adaptive (i.e., time-dependent) sampling of the spatial space ($x$) as well as the high-dimensional random space ($\xi$).  Overall, the computational cost of the proposed algorithm is $\mathcal{O}(r^2(n+s))$,  since $p \sim r$.  We refer to this algorithm as \emph{sparse TDB-ROM} (S-TDB-ROM) because the presented algorithm is not only to the DBO formulation and an identical procedure can be used to approximate $\bm{F}$ in DO and BO formulations. 

\subsection{Error Analysis}
In this section, we present error bounds for the low-rank approximation of the right hand side of the TDB equations. To this end, we first show that the above decomposition is equivalent to a CUR factorization of matrix $\bm{F}$. We then rely on existing CUR approximation error analyses to show that \cref{eq:F} is a near-optimal  approximation of $\bm{F}$ \cite{sorensen_deim_2016,mahoney_cur_2009}.  We present our analysis in the following Theorem.

\begin{theorem}\label{thm:CUR}
Let $\hat{\bm{F}} = \bm{U}_F\bm{Z}_F^T$ be a low-rank approximation of $\bm{F}$ according to the procedure explained in \S \ref{sec:sparse_interp}. Then  $\hat{\bm{F}}=\bm{U}_F\bm{Z}_F^T$ is equivalent to a CUR factorization of $\bm{F}$ such that the values of $\hat{\bm{F}}$ at all  sampled rows and columns are exact, i.e.,  $\hat{\bm{F}}(\bm{p},\bm{q}) = \bm{F}(\bm{p},\bm{q})$.
Proof. We first show that $\hat{\bm{F}} = \bm{U}_F\bm{Z}_F^T$ is equivalent to a CUR factorization of $\bm{F}$. To this end, it is sufficient to show that $\hat{\bm{F}}=\bm{C}\bm{U}\bm{R}$, where $\bm{C} = \bm{F}(:,\bm{q})$ and $\bm{R} = \bm{F}(\bm{p},:)$.  To show this, we replace $\bm{U}_F$ from \cref{eq:eig_F_q} and $\bm{Z}_F$ from \cref{eq:Z_f} into \cref{eq:F}, which results in:
\begin{equation*}
    \hat{\bm{F}} = \bm{F}(:,\bm{q})\bm{\Psi}_F\bm{\Lambda}_F^{-1/2} \bm{U}_F(\bm{p},:)^{-1}\bm{F}(\bm{p},:)
\end{equation*}
Now by letting $\bm{U} = \bm{\Psi}_F\bm{\Lambda}_F^{-1/2} \bm{U}_F(\bm{p},:)^{-1}$, we observe that $\hat{\bm{F}}=\bm{U}_F\bm{Z}_F^T = \bm{C}\bm{U}\bm{R}$. 

Now we show that the values of $\hat{\bm{F}}$ at all  sampled rows and columns are exact, i.e.,  $\hat{\bm{F}}(\bm{p},\bm{q}) = \bm{F}(\bm{p},\bm{q})$. To show this, from \cref{eq:F} we have:
    \begin{equation}\label{eq:F_pq}
        \hat{\bm{F}}(\bm{p},\bm{q}) = \bm{U}_F(\bm{p},:)\bm{Z}_F(\bm{q},:)^T.
    \end{equation}
    On the other hand, from \cref{eq:Z_f}, we obtain:
    \begin{equation}\label{eq:Z_fq}
        \bm{Z}_F(\bm{q},:)^T =  \bm{U}_F(\bm{p},:)^{-1}\bm{F}(\bm{p},\bm{q}).
    \end{equation}
     Replacing \cref{eq:Z_fq} into \cref{eq:F_pq} results in $\hat{\bm{F}}(\bm{p},\bm{q})=\bm{F}(\bm{p},\bm{q})$. This completes the proof. 
\end{theorem}

One of the implications of \cref{thm:CUR} is   that $\hat{\bm{F}}$ is  the \emph{oblique projection} of $\bm{F}$ onto the column space spanned by $\bm{U}_F$ and the row space spanned by $\bm{Z}_F$. These oblique projection operators are given by:
\begin{equation}\label{eq:PQ}
    \mathcal{P} = \bm{U}_F (\bm{P}^T \bm{U}_F)^{-1}\bm{P}^T \quad \mbox{and} \quad \mathcal{Q} = \bm{Q}(\bm{Z}_F^T\bm{Q} )^{-1}\bm{Z}_F^T.
\end{equation}
Now in the following we show that  $\hat{\bm{F}} = \mathcal{P} \bm{F} \mathcal{Q}$. Using the projection operators $\mathcal{P}$ and  $\mathcal{Q}$ from \cref{eq:PQ}, we have:
\begin{align*}
 \mathcal{P} \bm{F} \mathcal{Q}&=   \bm{U}_F (\bm{P}^T \bm{U}_F)^{-1}\bm{P}^T \bm{F} \bm{Q}(\bm{Z}_F^T\bm{Q} )^{-1}\bm{Z}_F^T\\
  &= \bm{U}_F \bm{U}_F(\bm{p},:)^{-1} \bm{F}(\bm{p},\bm{q}) \bm{Z}_F(\bm{q,:})^{-T} \bm{Z}_F^T. 
\end{align*}
From \cref{thm:CUR}, we have $\hat{\bm{F}}(\bm{p},\bm{q}) =\bm{F}(\bm{p},\bm{q}) $. Replacing $\hat{\bm{F}}(\bm{p},\bm{q})=\bm{U}_F(\bm{p},:)\bm{Z}^T_F(\bm{q},:)$   in the above equation results in 
\begin{align*}
 \mathcal{P} \bm{F} \mathcal{Q}=   \bm{U}_F  \bm{Z}_F^T = \hat{\bm{F}}.
\end{align*}
Now we use the results from \cite{sorensen_deim_2016}, where it was shown that for a CUR decomposition with the condition that $\hat{\bm{F}}(\bm{p},\bm{q}) = \bm{F}(\bm{p},\bm{q})$ the following holds:
\begin{equation}\label{eq:err_bnd}
     \|\bm{F}-\hat{\bm{F}} \|_2 \leq (\eta_p + \eta_q)\sigma_{p+1},  
\end{equation}
where $\sigma_{p+1}$ is the $(p+1)^{th}$ singular value of $\bm{F}$ and  $\eta_p =\|(\bm{P}^T\bm{U}_F)^{-1} \|_2$ and $\eta_q =\|(\bm{Q}^T\bm{Y}_F)^{-1} \|_2$.  The error bound given by \cref{eq:err_bnd} is valid irrespective of the choice of columns and rows ($\bm{p}$ and $\bm{q}$). Note that $\sigma_{p+1}$ is the optimal error, which is obtained if $\bm{F}$ is approximated with  a truncated rank-$p$ SVD of $\bm{F}$, i.e., $\|\bm{F}-\bm{F^*} \|_2 = \sigma_{p+1}$, where $\bm{F}^*$ is the optimal rank-$p$ approximation of $\bm{F}$. Therefore, $\eta_p+\eta_q$   is a scaling factor for the error and the goal of sparse a point selection algorithm is to minimize $\eta_p$ and $\eta_q$. Both DEIM and Q-DEIM point selection algorithms achieve near-optimal performance by ensuring that $\eta_p$ and $\eta_q$ remain small.  In this paper, we demonstrate the performance of both of these algorithms. We also note that there are other algorithms that can be used here. See for example algorithms based on leverage scores \cite{boutsidis_optimal_2014, drineas_relative-error_2008, mahoney_cur_2009}. 

\subsection{Rank-adaptive Approximation}\label{sec:adaptive}
To maintain the error of approximating $\bm{F}$ with $\hat{\bm{F}}$ below some desired threshold, it is natural to expect that the rank of the approximation ($p$) must change in time. Informed by the error analysis presented in the previous section, we propose an  algorithm for rank addition and removal. First, we note that the singular values of $\hat{\bm{F}}$ provide an approximation of the first $p$ singular values of $\bm{F}$. Since $\hat{\bm{F}}$ is a rank-$p$ approximation, we cannot compute $\sigma_{p+1}$ from \cref{eq:Y_f}. However, we know that $\sigma_p \geq \sigma_{p+1}$. Therefore, $\|\bm{F}-\hat{\bm{F}} \| \leq (\eta_p + \eta_q)\sigma_{p+1} \leq (\eta_p + \eta_q)\sigma_{p}$. As a result one can utilize the value of $\sigma_p$ as an indicator for the  approximation error. Since relative error is preferred to an absolute error, we propose to  use $\epsilon = \sigma_p^2/(\sum_{i=1}^p \sigma_i^2)$ as the criterion for rank addition/removal. We also use $\sigma_i^2$, because of its connection to the Frobenius norm for matrices, i.e., $\| \hat{\bm{F}} \|^2_F = \sum_{i=1}^p \sigma_i^2$. Setting a hard threshold may cause repetitive mode addition/removal. To prevent this unfavorable behavior, we propose to use \emph{buffer} interval for the error \cite{Babaee_Choi_Sapsis_Karniadakis_2017,Ashtiani_2022}. This means we set a lower bound ($\epsilon_l$) and an upper bound ($\epsilon_u$) for the error and  modes are added/removed to maintain $\epsilon_l \leq \epsilon \leq \epsilon_u$: If $\epsilon > \epsilon_u$,  we increase $p$ to $p+1$ and if $\epsilon < \epsilon_l$ we decrease $p$ to $p-1$. In the case of rank addition, we need to sample $p+1$ columns to compute $\bm{U}_F$. In  the DEIM and Q-DEIM algorithms   the number of rows or columns that can be computed is equal to the number of  singular vectors. However,  we use $\bm{Y}_F$ from the previous time step to sample the columns of $\bm{F}$, and $\bm{Y}_F$ has only $p$ singular vectors. To find the index for $(p+1)^{th}$ column,  we use the L-DEIM  algorithm \cite{LDEIM}, which  uses deterministic leverage scores -- facilitating  sampling more points   than the number of input singular vectors. The L-DEIM algorithm from \cite{LDEIM} is presented in \ref{AP_L-DEIM_Algorithm}.


\section{\label{sec:DC}Demonstration cases}

\subsection{Stochastic Burgers’ Equation}

\begin{figure}[!t]
    \centering
    \includegraphics[width=.7\textwidth]{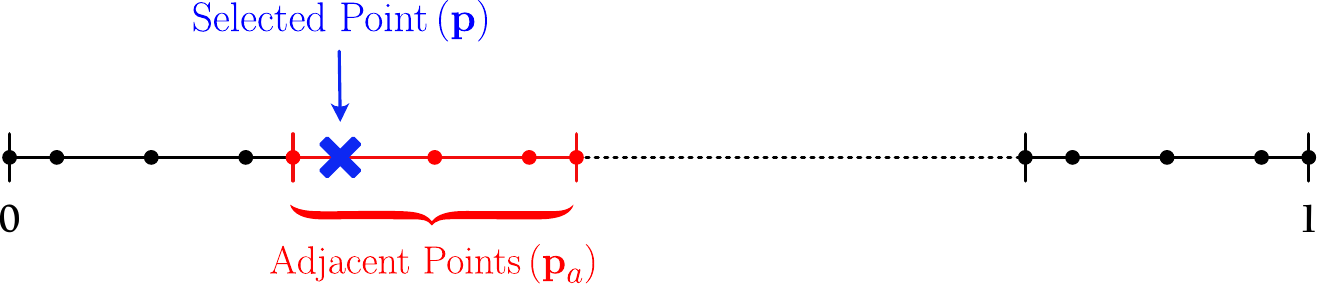}
    \caption{Burgers' equation: Schematic of the spectral element discretization and one of the selected points ($\bm{p}$) shown by a blue cross. One element is shown in red. The adjacent points ($\bm{p}_a$) required for computing the spatial derivatives at the selected point are shown by red circles. The Legendre-Gauss-Lobatto (LGL) scheme is used as collocation points. For this discretization, the solution at all collocation points in the element that selected point resides is needed. }
    \label{Burger_Schematic}
\end{figure}

For the first test case, we consider one-dimensional Burgers’ equation. As explained in \ref{AP_Computational Cost}, it is possible to achieve nominal TDB-ROM speedup for SPDEs with quadratic nonlinearity via an intrusive approach. In this demonstration, we show that S-TDB-ROM can achieve the nominal TDB-ROM speedup in a significantly less intrusive manner. We consider the Burgers' equation subject to random initial and boundary conditions as follows: 
\begin{equation*}
\begin{aligned}
\frac{\partial v}{\partial t}+v \frac{\partial v}{\partial x} &= \nu \frac{\partial^{2} v}{\partial x^{2}}, && x \in[0,1], t \in[0,1], \\
v(0, t; \xi) &= -\sin (2 \pi t)+\sigma_{t} \sum_{i=1}^{d} \lambda_{t_{i}}\varphi_{i}(t) \xi_{i}(\xi), && x=0, \xi_{i}\sim \mathcal{N}(0,1),\\
v(x, 0 ; \xi) &= 0.5\sin(2\pi x)\left(e^{\cos(2\pi x)}-1.5\right)+\sigma_{x} \sum_{i=1}^{d} \sqrt{\lambda_{x_{i}}} \psi_{i}(x) \xi_{i}(\xi), &&   x \in[0,1], \xi_{i}\sim \mathcal{N}(0,1),
\end{aligned}
\end{equation*}
where $\nu=0.05$ (except where otherwise stated) and we have Dirichlet boundary condition (BC) at $x = 1$ and stochastic Dirichlet BC at $x = 0$. The random space is taken to be $d = 4$ dimensional and $\xi_{i}$'s are sampled from a normal distribution with $s= 256$, $\sigma_{t}=0.01$, $\varphi_{i}(t)=\sin(i\pi t)$,  $\lambda_{t_{i}}=0.01/i^2$, and $\lambda_{x_{i}}$ and $\psi_{i}(x)$ are the eigenvalues and eigenvectors of the spatial squared-exponential kernel with $\sigma_{x}=0.005$, respectively. The fourth-order explicit Runge-Kutta method is used for time integration with $\Delta t=2.5 \times 10^{-4}$. For spatial discretization of the domain, the spectral element method and the Legendre-Gauss-Lobatto (LGL) scheme are used with 101 elements and polynomial order 4 which results in the total points in the domain being $n= 405$. \cref{Burger_Schematic} shows a schematic of the computational domain where the value of the adjacent points ($\bm{p}_a$) in the element are required for computing derivative at a selected point ($\bm{p}$), which is required for computing $\bm{Z}_F$. 

\begin{figure}[!t]
    \centering
    \includegraphics[width=.9\textwidth]{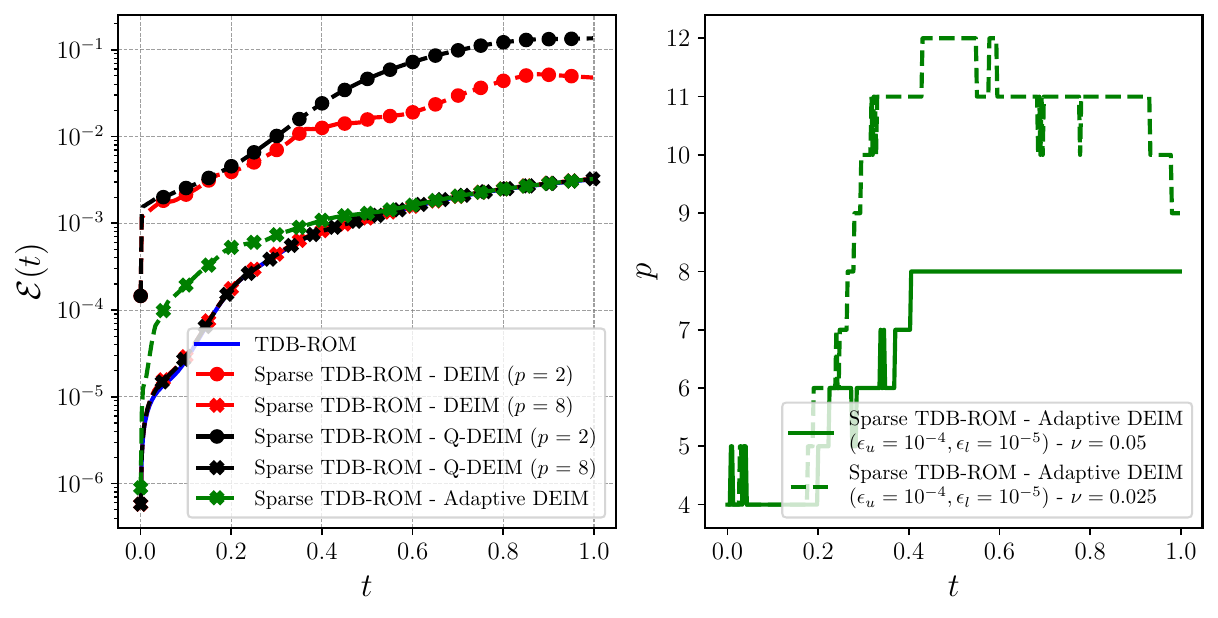}
    \caption{Stochastic Burgers' equation. Left: Error ($\mathcal{E}$) comparison for the TDB-ROM and S-TDB-ROM  for two different numbers of selected points $p = 2$ and $8$ and different sampling methods ($\nu=0.05$). Right: the number of  samples versus time for the adaptive DEIM method and two different values of $\nu$. The lower error bound ($\epsilon_l$) is equal to $10^{-5}$ and the upper error bound ($\epsilon_u$) is equal to $10^{-4}$ .}
    \label{fig:Error_Burger}
\end{figure}

The S-TDB-ROM and TDB-ROM solutions are compared with the FOM solution reduced to the same rank via KL decomposition. The total error is defined as:
\begin{equation*}
\mathcal{E}(t) =\left\|\mathbf{V}_{\text{TDB-ROM}}(t)-\mathbf{V}_{\text{FOM}}(t)\right\|_{F}.
\end{equation*}
The error between the TDB-ROM and S-TDB-ROM is due to the low-rank approximation of the right hand side of the SPDE. The total error i.e., $\mathcal{E}$, for the TDB-ROM and S-TDB-ROM methods with $r=5$ is shown in \cref{fig:Error_Burger}. As it can be seen, by increasing  $p$, the difference between the TDB-ROM and S-TDB-ROM becomes smaller, and with $p=8$ points, S-TDB-ROM is roughly as accurate as TDB-ROM. This indicates the presence of low-rank structure in the nonlinear term and the fact that the S-TDB-ROM method provides a good approximation for  $\mathbf{F}$. Also, both DEIM and Q-DEIM sampling strategies show similar accuracy.
On the other hand, for a lower error bound ($\epsilon_l$) equal to $10^{-5}$ and an upper error bound ($\epsilon_u$) equal to $10^{-4}$ in the adaptive DEIM method,  fewer number of points are selected   at the beginning of the simulation and as the system evolves, $p$ increases --- indicating the rank of the right hand side increases with time. Also, to show the effect of nonlinearity in S-TDB-ROM, we decreased $\nu$ from $\nu=0.05$ to $\nu=0.025$. It is evident from the right panel in \cref{fig:Error_Burger} that larger values of $p$ are required for the case of $\nu=0.025$. 

\begin{figure}[!t]
    \centering
    \includegraphics[width=\textwidth]{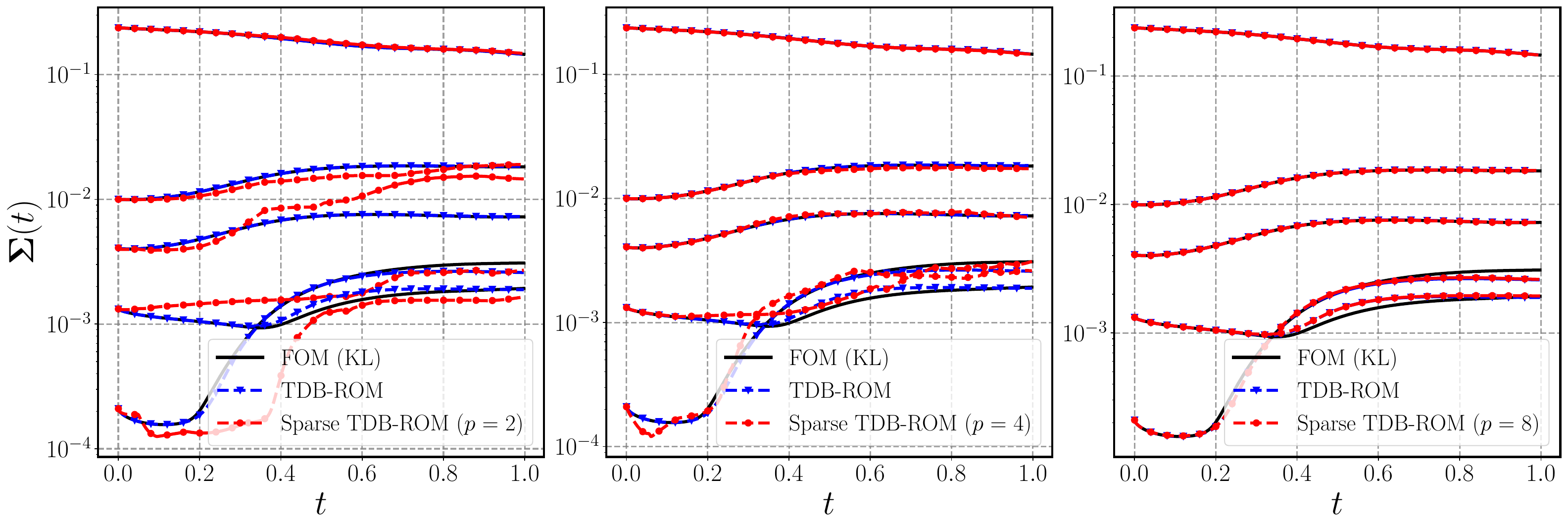}
    \caption{Stochastic Burgers' equation: Singular value comparison for FOM (KL), TDB-ROM, and S-TDB-ROM methods. The values are compared for three different numbers of selected points $p = 2, 4, \text{and}\, 8$ and the DEIM algorithm.}
    \label{fig:Sigma_Burger}
\end{figure}

In \cref{fig:Sigma_Burger} the instantaneous singular values of $\bm{\Sigma}(t)$ obtained from TDB-ROM, S-TDB-ROM and the $r$ largest singular values of the FOM solution (KL singular values) are shown.  The TDB-ROM solution shows a significant deviation from the FOM and TDB-ROM with $p=2$ which is due to the $\hat{\bm{F}}$  approximation error. However, as $p$ increases the singular values of TDB-ROM and S-TDB-ROM match. The deviation between the singular values of S-TDB-ROM and KL is due to the reduced order modeling error.

\begin{figure}[!t]
    \centering
    \includegraphics[width=\textwidth]{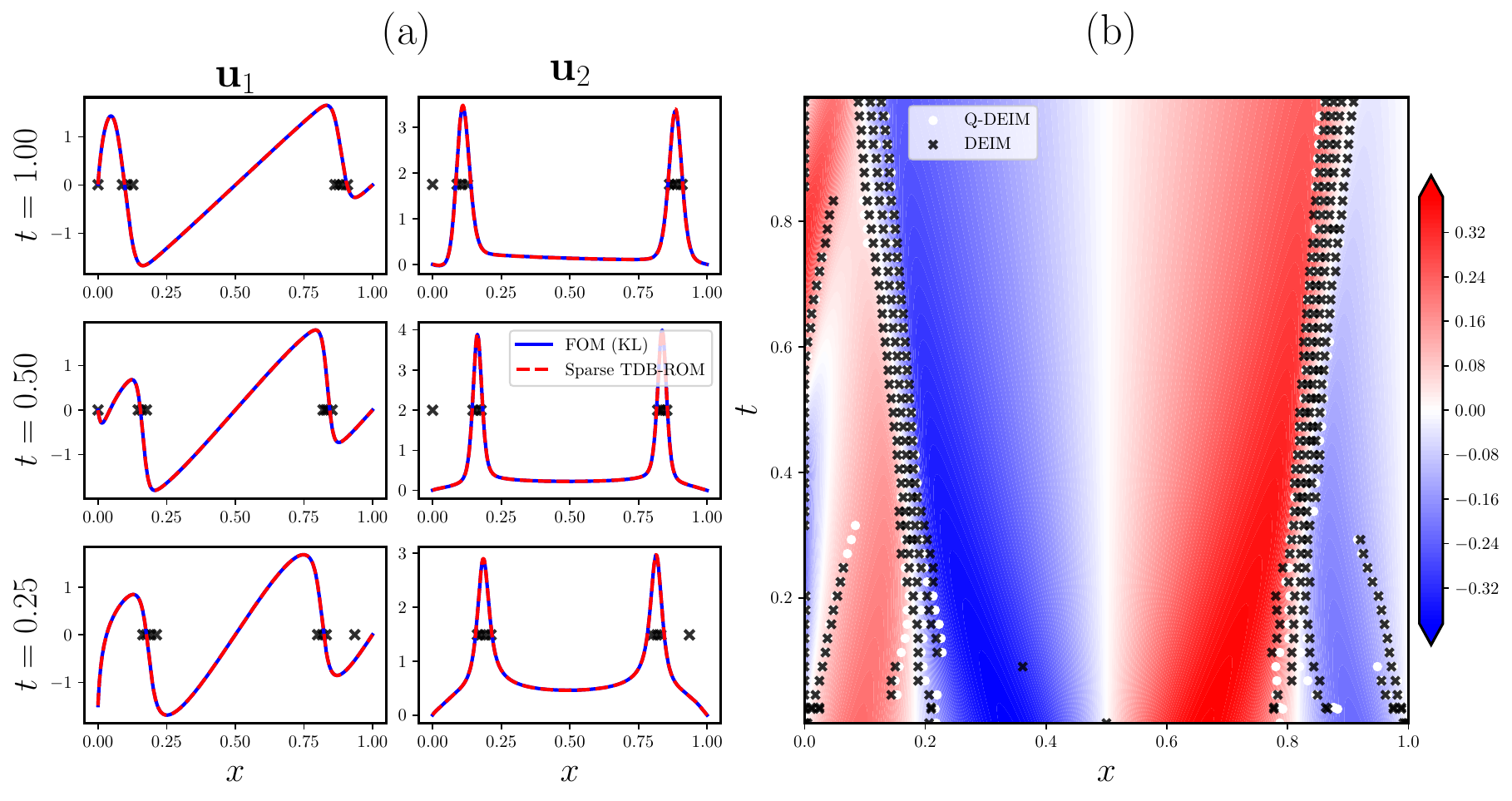}
    \caption{Stochastic Burgers' equation: (a) Evolution of the S-TDB-ROM and FOM (KL) first two spatial modes ($\mathbf{u}_{1}$, $\mathbf{u}_{2}$) at $t=0.25,0.50,1.00$. Crosses show    the selected points by the DEIM algorithm, and (b) contour plot of the S-TDB-ROM solution with the distribution of the selected spatial points usins DEIM and Q-DEIM algorithms with $p = 8$.}
    \label{fig:Points_Burger}
\end{figure}

\cref{fig:Points_Burger} shows the evolution of the two most dominant modes of the S-TDB-ROM and FOM (KL). These modes are ranked based on the instantaneous singular values, which implies that $\mathbf{u}_{1}$ and  $\mathbf{u}_{2}$ are the two most energetic modes. Excellent  agreement between S-TDB-ROM and FOM (KL) is observed. Also, we can see the distribution of the selected spatial points in the contour plot. The sparse points selected by DEIM and Q-DEIM algorithms in the physical space are very similar to each other and they concentrate near the two shocks where  there is a high gradient in the solution. Because we use a TDB method here, the location of selected points varies in each time-step. The left boundary is stochastic and interestingly the boundary point is not always selected. At the right boundary deterministic Dirichlet boundary condition is imposed and this boundary is only selected in the few first time steps.

\begin{figure}[!t]
   \centering
   \includegraphics[width=\textwidth]{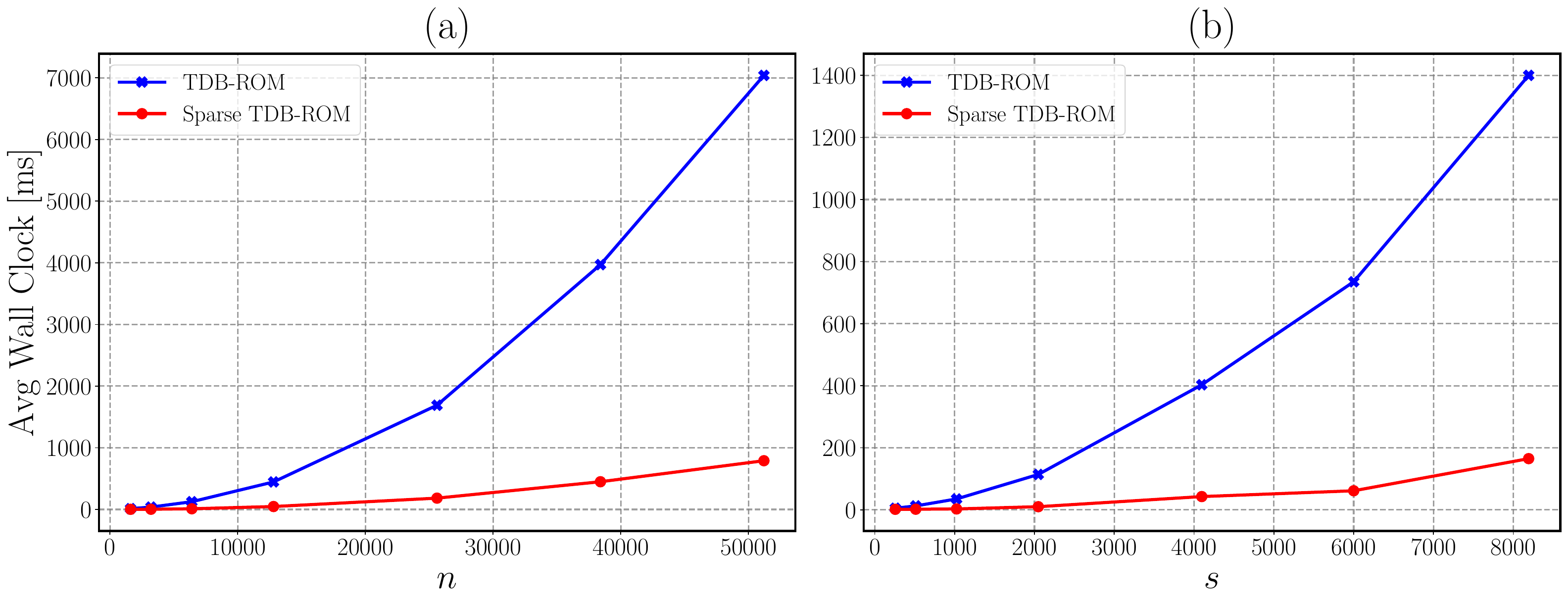}
   \caption{Stochastic Burgers' equation: Average wall clock comparison for different (a) the number of points ($n$) and (b) number of samples ($s$) between TDB-ROM and S-TDB-ROM with the DEIM algorithm.}
   \label{fig:ne_ns_Burger}
\end{figure}

 \cref{fig:ne_ns_Burger} compares the CPU time of the S-TDB-ROM versus TDB-ROM. For the S-TDB-ROM, we fix $p=8$, $r=5$, and $s=256$ variables in (a) and $p=8$, $r=5$, and $n=405$ variables in (b).
 Clearly, as $n$ or $s$ becomes larger, the disparity between CPU time of S-TDB-ROM and TDB-ROM increases.

\subsection{Stochastic Compressible Navier-Stokes Equations}
In the second demonstration, we apply S-TDB-ROM to compressible Navier-Stokes Equations, which has  non-polynomial nonlinearities. As a result, solving the TDB-ROM using \cref{eq:DBO_evol_S,eq:DBO_evol_U,eq:DBO_evol_Y} is as expensive as that of solving FOM. We consider two cases here. In the first case, we consider a small value of $s$ so that we can compare S-TDB-ROM with TDB-ROM and FOM. In the second case, we solve the compressible flow subject to 100-dimensional random perturbations.  For this case, we consider a large value of $s$, for which we could not run  FOM nor TDB-ROM given the computational resources at our disposal. But we show that we could solve the same problem using S-TDB-ROM methodology on an NVIDIA QUADRO P5000 GPU card with 2560 CUDA cores and 16 GB memory. The 2D compressible Navier-Stokes equations are given by:
\begin{equation*}
\begin{aligned}
\frac{\partial \rho}{\partial t}+\frac{\partial \rho v_{j}}{\partial x_{j}} &=0, \\
\frac{\partial \rho v_{i}}{\partial t}+\frac{\partial \rho v_{i} v_{j}}{\partial x_{j}} &=-\frac{\partial p}{\partial x_{i}}+\frac{\partial \tau_{i j}}{\partial x_{j}}, \\
\frac{\partial E}{\partial t}+\frac{\partial\left(E v_{j}\right)}{\partial x_{j}} &=-\frac{\partial p v_{j}}{\partial x_{j}}+\frac{\partial\left(\tau_{i j} v_{i}\right)}{\partial x_{j}}-\frac{\partial q_{j}}{\partial x_{j}},
\end{aligned}
\end{equation*}
where the temperature $T(x, t)$, pressure $p(x, t)$, velocity $v(x, t)$, total energy $E(x, t)$, and density $\rho(x, t)$ are primary transport variables. The viscosity flux $\tau$, and heat flux $q$ are defined as:
\begin{equation*}
\tau_{i j} =\frac{1}{\operatorname{Re}}\left(\frac{\partial v_{i}}{\partial x_{j}}+\frac{\partial v_{j}}{\partial x_{i}}-\frac{2}{3} \frac{\partial v_{k}}{\partial x_{k}} \delta_{i j}\right)\quad \mbox{and} \quad 
q_{j} =-\frac{1}{Ec \cdot Pe} \frac{\partial T}{\partial x_{j}},
\end{equation*}
where $e$ is the internal energy, $E= \rho e+\frac{1}{2} \rho v_{i} v_{i}$ is the total energy, $Pe=\operatorname{Re} . Pr$, $Ec=(\gamma-1) Ma^{2}$, and $Ma$ are Peclet, Eckert, and Mach numbers, respectively. In our simulations, $\operatorname{Re}=3000$, $Pr=1.0$, $\gamma=1.4$, and $Ma=0.5$. The fourth-order explicit Runge-Kutta method is utilized for time integration with $\Delta t=5 \times 10^{-4}$. Periodic boundary conditions are imposed on all four boundaries and initial pressure on the entire domain is set to $p(x,y,0)=1$. Also, initial temperature and velocity are defined as follows:
\begin{equation*}
\begin{aligned}
        u(x, y, 0) &= \overline{u} +  \frac{2L_{x}\delta}{h^{2}} ((y-b)e^{\frac{-(y-b)^2}{h^{2}}}+(y-a)e^{\frac{-(y-a)^2}{h^{2}}})\sin(\frac{10\pi x}{L_{x}}),\\
        v(x, y, 0) &= 10\pi\delta (e^{\frac{-(y-b)^2}{h^{2}}}+e^{\frac{-(y-a)^2}{h^{2}}})\cos(\frac{10\pi x}{L_{x}}),\\
        T(x, y, 0) &= 0.5+0.25(\tanh(\frac{y-y_{min}}{h}) - \tanh(\frac{y-y_{max}}{h})),
\end{aligned}
\end{equation*}
where $\delta=\num{1.45d-4}$, $y_{min}=0.45$, $y_{max}=0.55$, $U_{max}=1$, $a=0.45$, $b=0.55$, $h=0.01$, and
\begin{equation*}
\overline{u} = 0.5U_{max}(\tanh(\frac{y-y_{min}}{h})-\tanh(\frac{y-y_{max}}{h})-1).
\end{equation*}
We solve the flow with the above initial condition for three units of time ($t_s=3$). Then, the TDB-ROM and S-TDB-ROM  are initialized with the transport variables at $t_{s} = 3$ as well as random fluctuations as shown below:
\begin{equation*}
\begin{aligned}
u(x, y, t_{s};\xi)&= u(x, y, t_{s}) + \sum_{k=1}^{d} \lambda_{{k}} ((y-b)e^{\frac{-(y-b)^2}{h^{2}}}+(y-a)e^{\frac{-(y-a)^2}{h^{2}}})\sin(\frac{2k\pi x}{L_{x}}) \xi_{k}(\xi),\\
v(x, y, t_{s}; \xi)&= v(x, y, t_{s}) + \sum_{k=1}^{d} \lambda_{{k}} ((y-b)e^{\frac{-(y-b)^2}{h^{2}}}+(y-a)e^{\frac{-(y-a)^2}{h^{2}}})\cos(\frac{2k\pi x}{L_{x}}) \xi_{k}(\xi),
\end{aligned}
\end{equation*}
where $\lambda_{{k}}=10/k^2$ and $\xi_{k}$ are independent random variables sampled from a normal distribution.

The domain is discretized using the finite difference method on a uniform $256\times 256$ grid where $L_x = 2$ and $L_y = 1$. Also, the four  variables $\rho$, $\rho u$, $\rho v$, and $E$ are stacked together  in order to create a global mode, i.e., $\bm{v} =[\rho, \  \rho u, \  \rho v, \ E ]^T  \in \mathbb{R}^{n}$, where  $n = 4\times 256^2=262144$.

\begin{figure}[!t]
    \centering
    \includegraphics[width=.65\textwidth]{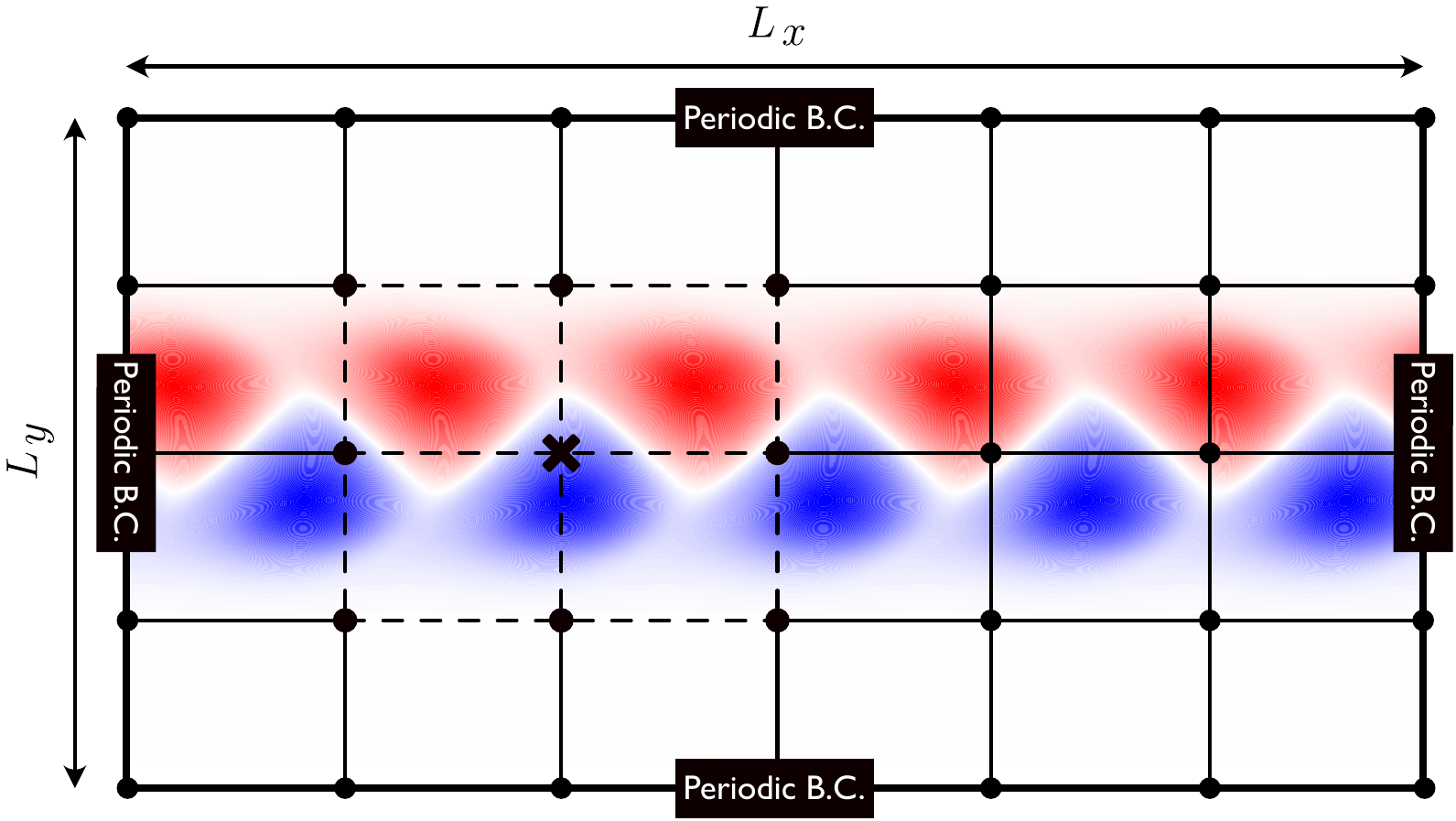}
    \caption{Stochastic compressible Navier-Stokes equations: Schematic of the domain, and required adjacent points $\bf{p_a}$ (dots on the dash lines) for computing derivative at a selected point $\bf{p}$ (cross) in the finite difference method.} 
    \label{fig:NS_Schematic}
\end{figure}

As mentioned before, to sample rows of $\bm{F}$, we need to calculate the spatial derivative at the selected  points. This requires the values of $\bm{V}$ at the adjacent points.  In  \cref{fig:NS_Schematic}  the schematic of the adjacent points (shown by filled circles) required for computing the derivative at a selected point (shown by a cross) is shown. For each selected point, values of the solution at 8 adjacent points must be provided.

\begin{figure}[!t]
    \centering
    \includegraphics[width=.9\textwidth]{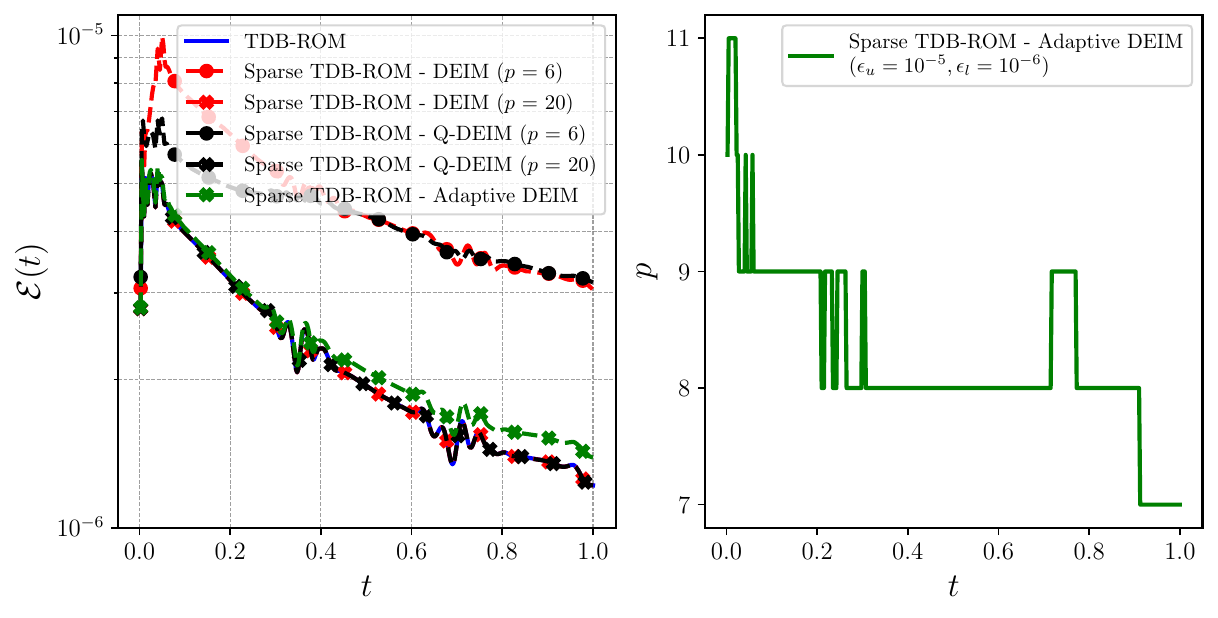}
    \caption{Stochastic compressible Navier-Stokes equations: Error ($\mathcal{E}$) comparison for the TDB-ROM and S-TDB-ROM ($d=20$, $s=150$, and $r=5$) as compared with the FOM solution for two different numbers of selected points $p = 6$ and $20$  and different sampling methods. The lower error bound ($\epsilon_l$) is equal to $10^{-6}$ and the upper error bound ($\epsilon_u$) is equal to $10^{-5}$ for the adaptive DEIM method.}
    \label{fig:Error_NS}
\end{figure}

In the first case that we present, $d=20$ and we consider $r=5$ modes.  We choose a grossly insufficient number of samples  $s=150$. However, our goal here is to assess the performance of S-TDB-ROM, which requires solving 2D compressible Navier-Stokes equations for 150 samples as well as solving the TDB-ROM equations.  In \cref{fig:Error_NS}, we depict the total error of the solution versus time for both the TDB-ROM and S-TDB-ROM. We use different numbers of samples ($p$). As illustrated in the previous test case, the S-TDB-ROM is able to obtain highly accurate results, comparable with the TDB-ROM method, with a few selected points at a significantly reduced cost.  Also, while we cannot observe any noteworthy difference between the DEIM and Q-DEIM methods, for a lower error bound ($\epsilon_l$) equal to $10^{-6}$ and an upper error bound ($\epsilon_u$) equal to $10^{-5}$, the adaptive DEIM method approximately achieves the same level of accuracy with a considerably smaller number of points compared to the fixed 20 points used in the DEIM and Q-DEIM algorithms.

\begin{figure}[!t]
    \centering
    \includegraphics[width=\textwidth]{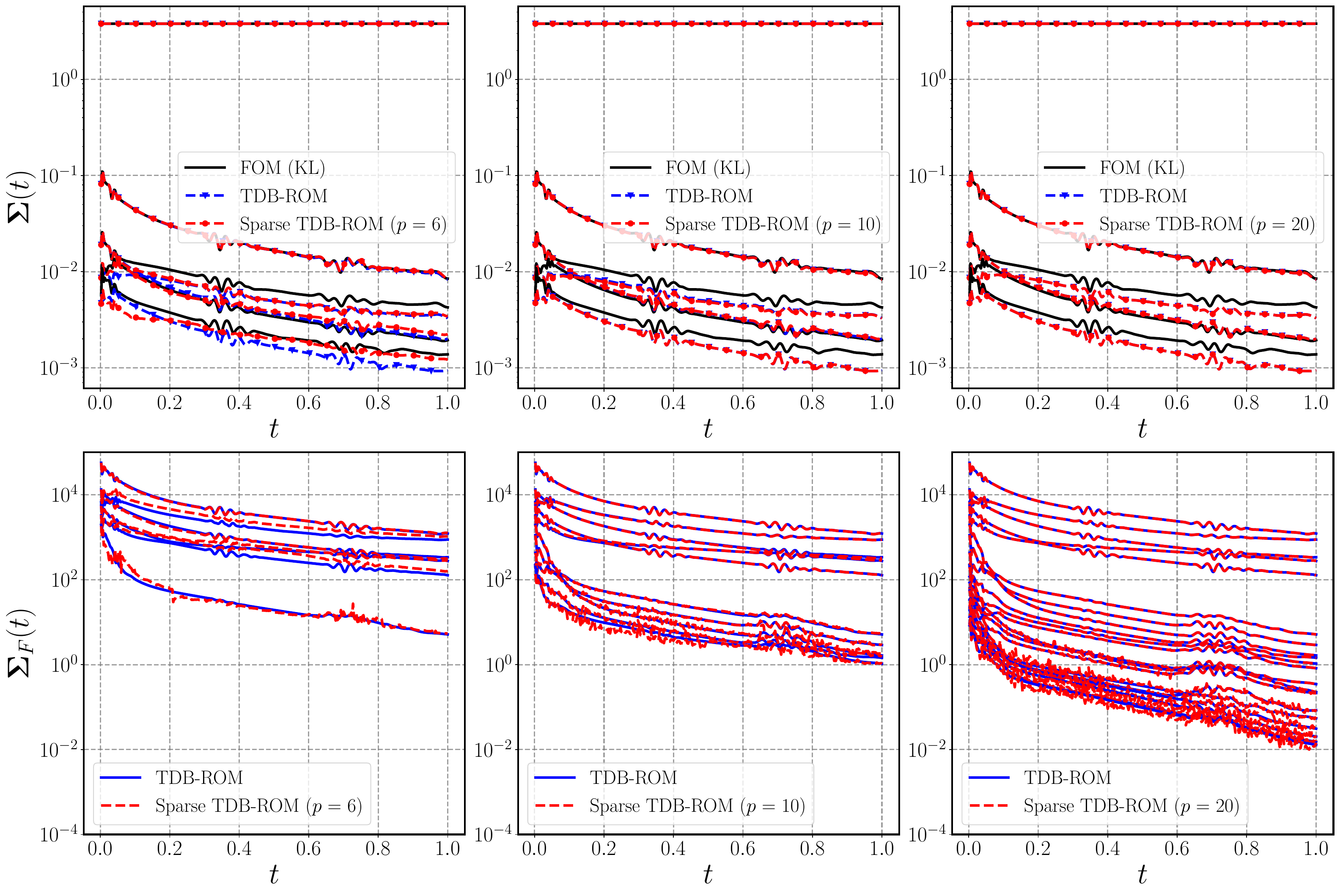}
    \caption{Stochastic compressible Navier-Stokes equations: Comparison of the singular values and nonlinear term ($\mathbf{F}$) singular values for the FOM, TDB-ROM, and S-TDB-ROM methods ($d=20$, $s=150$, and $r=5$). The values are compared for three different numbers of selected points of $p = 6, 10, \text{and}\, 20$.}
    \label{fig:Sigma_NS}
\end{figure}

In the first row of \cref{fig:Sigma_NS} instantaneous singular values of the TDB-ROM and S-TDB-ROM methods for $p = 6$, 10, and 20 are shown. As $p$ increases the lower singular values of  S-TDB-ROM also match with those  of the TDB-ROM. In the second row of \cref{fig:Sigma_NS},  instantaneous singular values  of $\mathbf{F}$ in the TDB-ROM are compared against  $\hat{\bm{F}}=\bm{U}_F \bm{Z}_F^T$ from the S-TDB-ROM which represents the approximation of $\mathbf{F}$. We note that  the rank of $\bm{F}$ is more than the rank of $\hat{\bm{V}}$, which is $r=5$. This is because $\mathcal{F}$ is a nonlinear map, and the rank of $\bm{F}=\mathcal{F}(\hat{\bm{V}})$ is not in general the same as rank of $\hat{\bm{V}}$. The singular values of S-TDB-ROM  closely follow those of  TDB-ROM up to a certain $p$ and the reminder singular values generally show noise-like behavior. 

\begin{figure}[!t]
    \centering
    \includegraphics[width=\textwidth]{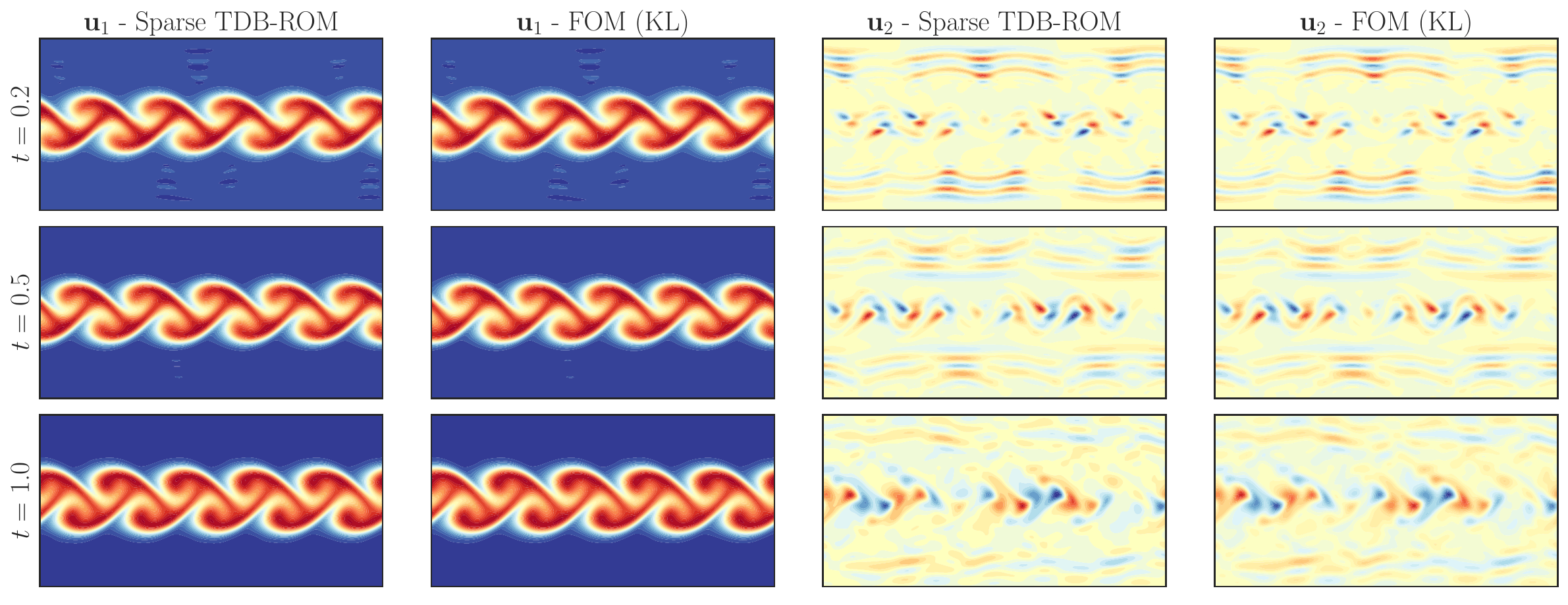}
    \caption{Stochastic compressible Navier-Stokes equations: Evolution of the first two spatial modes ($\mathbf{u}_{1}$, $\mathbf{u}_{2}$) of the S-TDB-ROM and FOM (KL) at $t=0.2,0.5,1.0$.}
    \label{fig:modes}
\end{figure}

In \cref{fig:modes}, the first two dominant spatial modes of the S-TDB-ROM and FOM (KL) at different time instants are shown. The KL modes  are computed by taking SVD of the FOM solution (for $s=150$ samples) at the time instance in question. It is clear that there is a good match between the S-TDB-ROM and FOM (KL) spatial modes.

\begin{figure}[!t]
    \centering
    \includegraphics[width=\textwidth]{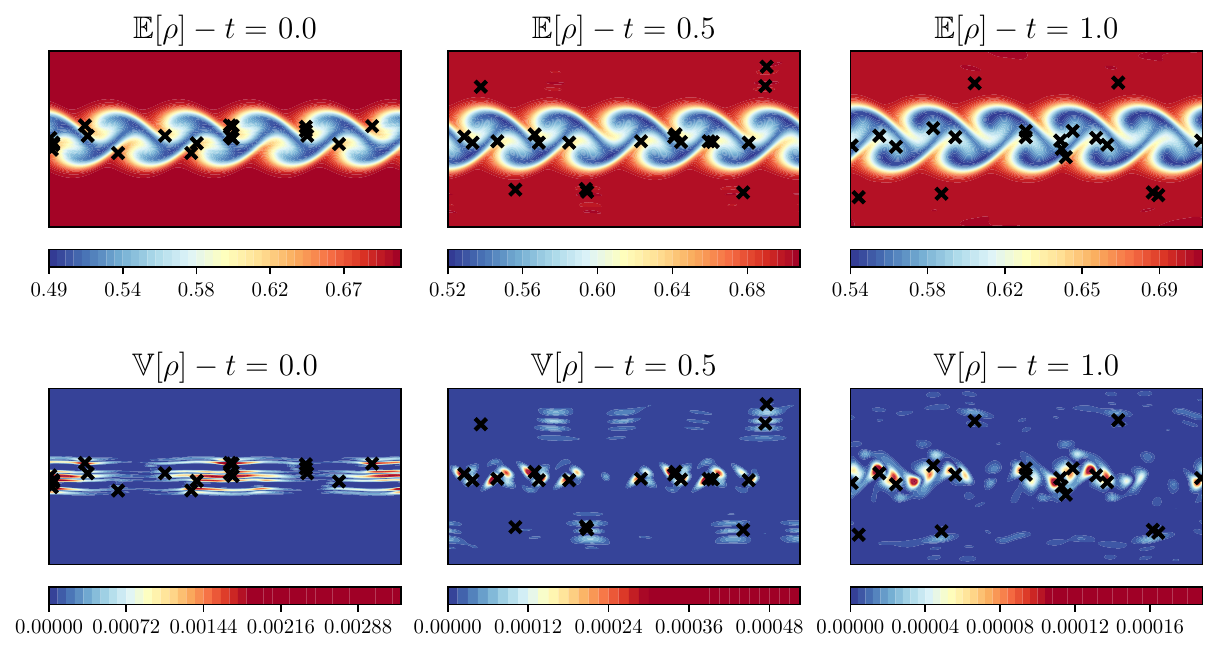}
    \caption{Stochastic compressible Navier-Stokes equations: Mean and variance of the density in different times and the DEIM selected points ($n = 262144$, $s=100000$, $d=100$, $r=7$, and $p=20$).}
    \label{fig:mean_var}
\end{figure}

In the second case, we consider $d=100$-dimensional random space and we choose $s=100000$ samples. We consider  $r=7$, and $p=20$. This case is particularly of interest since it shows the true capability of the S-TDB-ROM  for a case where  TDB-ROM  is prohibitively expensive to run. For running this setting with the TDB-ROM we would have to have available memory for storing a $262144\times 100000$ matrix ($\bm{F}$)  and compute the nonlinear map of a matrix of this size each time step ($\bm{F}=\mathcal{F}(\hat{\bm{V}})$). However, using the S-TDB-ROM method, we never need to form a matrix larger than $262144\times 20$. \cref{fig:mean_var} shows the mean and variance of the density at different times with selected points. The DEIM algorithm  selects  points in the regions with a high variance. This demonstrates the effectiveness of -- and the need for --  adaptive sampling where the DEIM samples  change in time according to the state of the solution. 

\begin{figure}[!t]
    \centering
    \includegraphics[width=\textwidth]{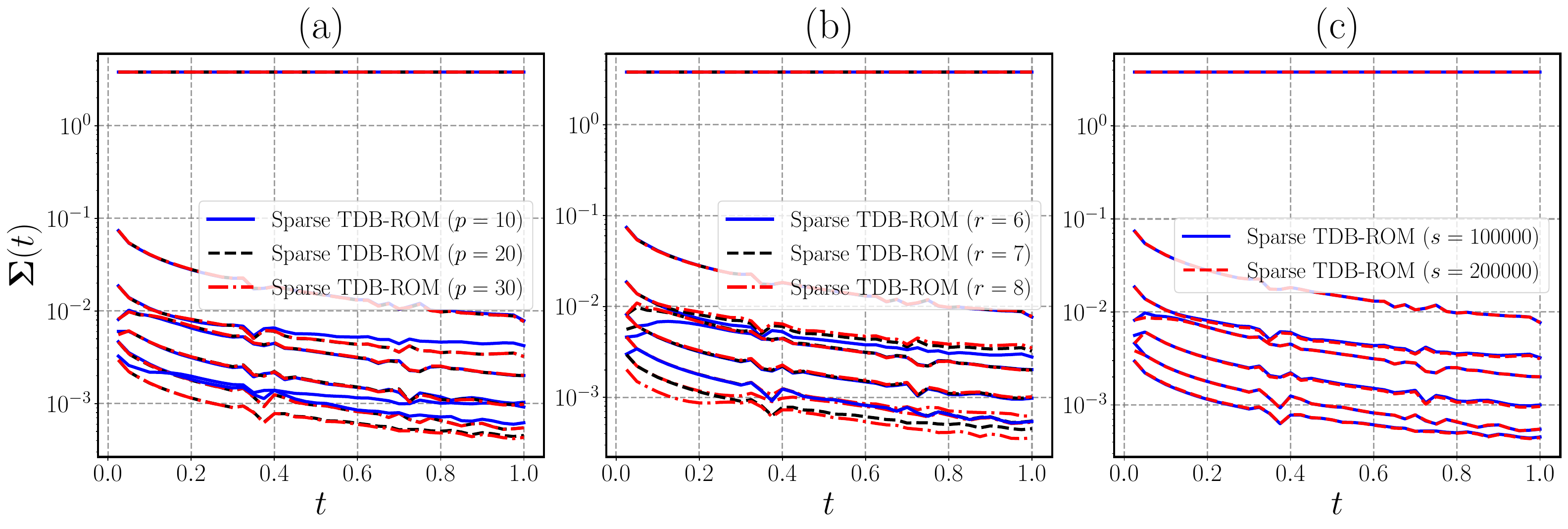}
    \caption{Stochastic compressible Navier-Stokes equations: Convergence of the S-TDB-ROM singular values for (a) different numbers of selected points $p = 10, 20, \text{and}\, 30$, (b) different numbers of modes $r = 6, 7, \text{and}\, 8$, and (c) different numbers of samples $s = 100000, \text{and}\, 200000$.}
    \label{fig:Sigma_NS_100000}
\end{figure}

\begin{figure}[!t]
    \centering
    \includegraphics[width=\textwidth]{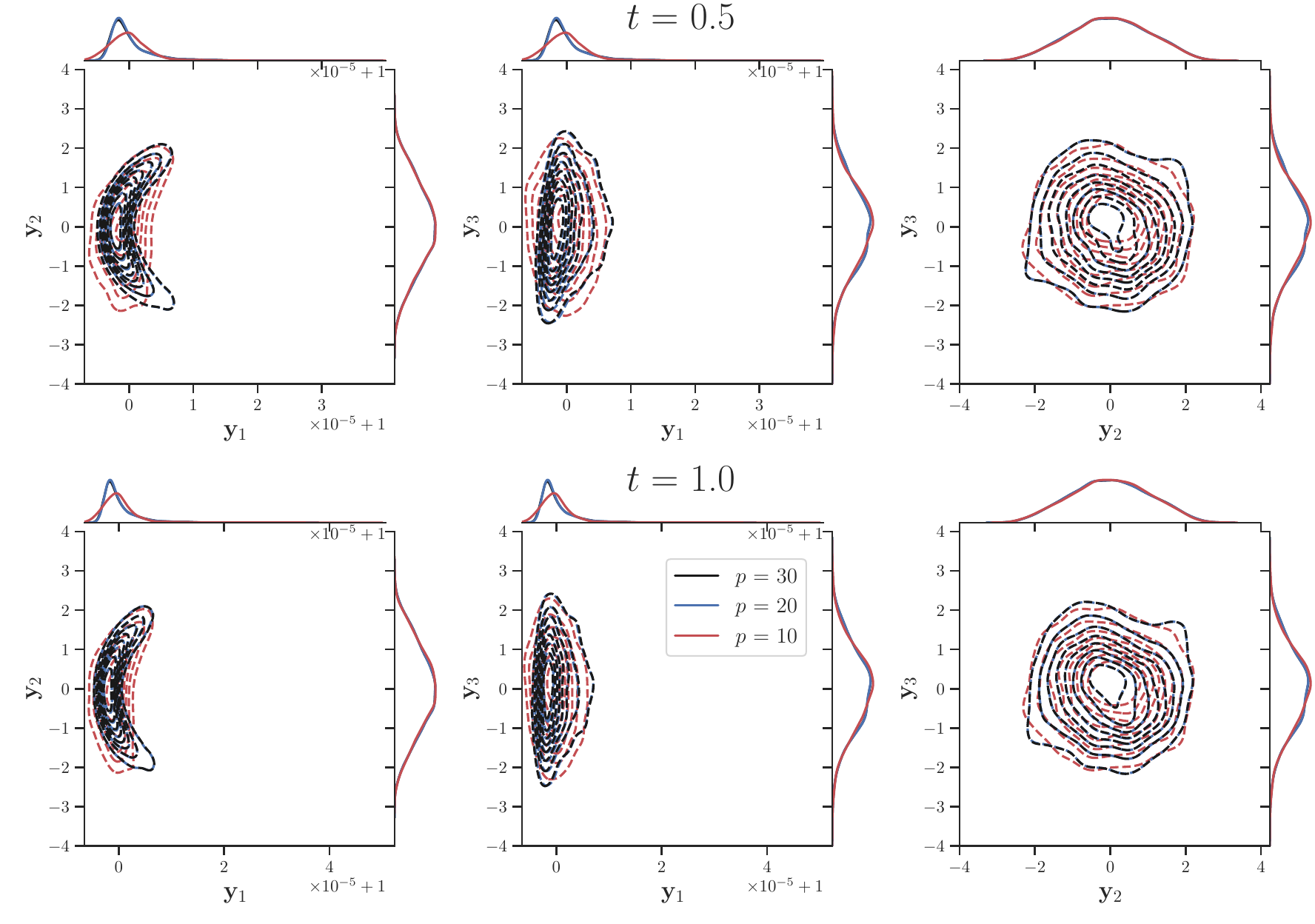}
    \caption{Stochastic compressible Navier-Stokes equations: Joint and marginal pdfs for modes 1, 2, and 3 of $\mathbf{Y}$. The values are compared for three different numbers of selected points $p = 10, 20, \text{and}\, 30$.}
    \label{fig:Yjoint_np}
\end{figure}

\begin{figure}[!t]
    \centering
    \includegraphics[width=\textwidth]{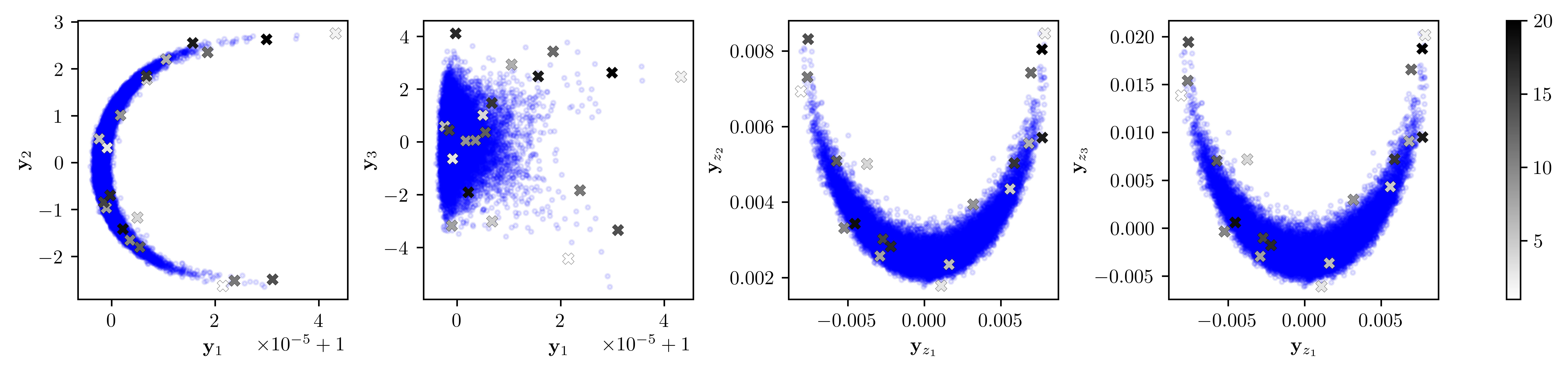}
    \caption{Stochastic compressible Navier-Stokes equations: Scatter plots for modes 1, 2, and 3 of $\mathbf{Y}$ and $\mathbf{Y}_{z}$ ($p = 20$). The color bar shows the index number of the selected points. The points with lower index numbers are  selected first.}
    \label{fig:Y_Scatter}
\end{figure}

The accuracy of the S-TDB-ROM method has been confirmed in previous cases by error comparison with the TDB-ROM method. However, we cannot use the TDB-ROM method in this case since the computational cost is prohibitive. Here we perform a convergence study by changing  $r$, $p$, and $s$, which are the reduced-order modeling parameters ($r,p$) as well as the number of samples $(s)$.  In each case, the two other parameters are fixed ($r=7$, and $s=100000$ for the convergence study of $p$, $r=7$, and $p=20$  for the convergence study of $s$, and $p=20$, and $s=100000$ for convergence study of $r$). From \cref{fig:Sigma_NS_100000} we have converged singular values with $r=7$, $p=20$, and $s=100000$. Furthermore, \cref{fig:Yjoint_np} depicts the convergence of the joint and marginal probability density functions (pdf)   for the first three dominant stochastic modes  ($\mathbf{y}_1$, $\mathbf{y}_2$, and $\mathbf{y}_3$). From \cref{fig:Yjoint_np} and \cref{fig:Y_Scatter}, the nonlinear relation between the first three $\mathbf{y}$ modes can be observed where the S-TDB-ROM method selects  points in both high and low probability with a specific order (shown with the color bar) to be able to interpolate the $\mathbf{Y}$ matrix with high accuracy.

\section{\label{sec:Conclusion}Conclusion}

We present a  methodology to reduce the computational cost of evaluating the right hand side term from $\mathcal{O}(ns)$ to $\mathcal{O}(r^2(n+s))$ for both non-homogeneous linear equations as well as any nonlinear SPDEs with generic nonlinearity (polynomial or non-polynomial). Moreover, the presented approach replaces the highly intrusive steps often done in linear and quadratic SPDEs with a procedure that is agnostic to the  type of the equation and in that sense it significantly reduces the level of intrusiveness of the derivation and implementation of the TDB evolution equations for different SPDEs. The  algorithm presents  a DEIM-based sparse interpolation strategy for the rows and columns of the right hand side of the SPDE. However, unlike the DEIM algorithm, the presented methodology does not require an offline data-driven step for the computation of the POD bases for the nonlinear terms. Doing so would detract from some of the key advantages of reduced order modeling based on TDBs. Also, we proposed a procedure for adaptively selecting points in different time steps by rank addition and removal according to a specified threshold for the error. This allows the algorithm to choose the number of required sampling points by their significance at different times.

We demonstrated the performance of the presented method on  two different case studies; the stochastic Burgers’ equation and the stochastic compressible Navier-Stokes equations. For small number of samples we showed that sparse TDB-ROM and the TDB-ROM in the decompressed form yield similarly accurate results for a large enough interpolation points. For the compressible Navier-Stokes equation, we considered a case with $10^5$ samples for which we could not solve the FOM not the TDB-ROM equations in the decompressed form using the computational resources at our disposal. However, we showed that we can solve this problem using the presented algorithm. 

\section*{\label{sec:DCI}Declaration of competing interest}
The authors declare that they have no known competing financial interests or personal relationships that could have appeared to influence the work reported in this paper.

\section*{\label{sec:Acknowledgement}Acknowledgement}
This work is sponsored by   the Air Force Office of Scientific
Research award (PM: Dr. Fariba Fahroo) FA9550-21-1-0247 and the National Science Foundation (NSF), USA under Grant CBET2042918. 

\appendix

\section{\label{AP_Sparse_Sampling}Sparse Sampling Methods}
\subsection{\label{AP_DEIM_Algorithm}Direct Empirical Interpolation Method (DEIM)}

The DEIM method seeks to find near optimal interpolation points for approximating a function versus a set of orthonormal bases ($\bm{\Psi}_p$). The DEIM pseudocode is presented via \cref{alg:DEIM}  and we refer to \cite{chaturantabut_nonlinear_2010} for more details on the DEIM framework.

\begin{algorithm}
\SetAlgoLined
\KwIn{$\mathbf{\Psi}_{p}=\left[\begin{array}{llll}\mathbf{\psi}_{1} & \mathbf{\psi}_{2} & \cdots & \mathbf{\psi}_{p}\end{array}\right]$}
\KwOut{$I_{p}$}
$\left[\rho, I_{1}\right]=\max \left|\mathbf{\psi}_{1}\right|$ \hspace{1.3cm} $\rhd$ choose the first index\;
$\mathbf{P}_{1}=\left[\mathbf{e}_{I_{1}}\right]$ \hspace{2.45cm} $\rhd$ construct first measurement matrix\;
\For{$i=2$ \KwTo $p$}{
$\mathbf{P}_{i}^{T} \mathbf{\Psi}_{i} \mathbf{c}_{i}=\mathbf{P}_{i}^{T} \mathbf{\psi}_{i+1}$ \hspace{0.46cm}  $\rhd$ calculate $c_{i}$\;
$\mathbf{R}_{i+1}=\mathbf{\psi}_{i+1}-\mathbf{\Psi}_{i} \mathbf{c}_{i}$ \hspace{0.24cm}  $\rhd$ compute residual\;
$\left[\rho, I_{i}\right]=\max \left|\mathbf{R}_{i+1}\right|$ \hspace{0.38cm} $\rhd$ find index of maximum residual\;
$\mathbf{P}_{i+1}=\left[\begin{array}{ll}\mathbf{P}_{i} & \mathbf{e}_{I_{i}}\end{array}\right]$ \hspace{0.38cm} $\rhd$ add new column to measurement matrix\;}
\caption{DEIM Algorithm \cite{chaturantabut_nonlinear_2010}}
\label{alg:DEIM}
\end{algorithm}

\subsection{\label{AP_Q-DEIM_Algorithm} Q-DEIM Algorithm}

While the DEIM algorithm is an efficient method for approximation of a nonlinear function, there are other approaches that are equally efficient. Q-DEIM was presented  \cite{drmac_new_2016} as a new method for for selecting interpolation points using the QR factorization with column pivoting.  This method has been established as a robust alternative framework for sensor placement in many applications \cite{data-driven_2018,chellappa_training_2021,peherstorfer_stability_2020}. The availability of the pivoted QR implementation in many open-source packages  makes this algorithm an efficient alternative for sparse sampling. \cref{alg:QDEIM} can replace the DEIM algorithm to construct the $I_{p}$.

\begin{algorithm}
\SetAlgoLined
\KwIn{$\mathbf{\Psi}_{p}=\left[\begin{array}{llll}\mathbf{\psi}_{1} & \mathbf{\psi}_{2} & \cdots & \mathbf{\psi}_{p}\end{array}\right]$}
\KwOut{$I_{p}$}
$\left[\mathrm{q}, \mathrm{r}, \mathrm{pivot}\right] \leftarrow$ \texttt{qr($\mathbf{\Psi}^T_{p}$)} \hspace{0.45cm}  $\rhd$ QR factorization with column pivoting\;
$I_{p}=\mathrm{pivot}(1:\mathrm{p})$ \hspace{1.49cm} $\rhd$ selecting first p elements of the pivot\;
\caption{Q-DEIM Algorithm \cite{drmac_new_2016}}
\label{alg:QDEIM}
\end{algorithm}

\subsection{\label{AP_L-DEIM_Algorithm} L-DEIM Algorithm}

One of the limitations of the DEIM algorithm is that the number of column indices that can be computed is restricted to the input singular vectors. Using L-DEIM, the number of  sample points could be larger than the number of  singular vectors ($p$). \cref{alg:LDEIM} shows the pseudocode of the proposed method and more details can be found in \cite{LDEIM}.
\begin{algorithm}
\SetAlgoLined
\KwIn{$\mathbf{\Psi} \in \mathbb{R}^{m \times p}$, target rank = $\hat{p}$}
\KwOut{$I_{\hat{p}}$}
\For{$i=1$ \KwTo $p$}{
$I_{i}=\operatorname{argmax}_{1 \leq j \leq m}\left|(\mathbf{\Psi}_{i})_{j}\right|$\;
$\mathbf{\Psi}_{i}=\mathbf{\Psi}_{i}-\mathbf{\Psi}(:,1:i) \cdot(\mathbf{\Psi}(I,1:i) \backslash \mathbf{\Psi}(I,i+1))$\;}
\For{$i=1$ \KwTo $m$}{
$\ell_{i}=\left\|[\mathbf{\Psi}]_{i:}\right\|$}
$\operatorname{sort}(\ell)$\;
Remove entries in $\ell$ corresponding to the indices in $I$\;
$I^{\prime}=\hat{p}-p$ indices corresponding to $\hat{p}-p$ largest entries of $\ell$\;
$I_{\hat{p}}=\left[I; I^{\prime}\right]$\;
\caption{L-DEIM Algorithm \cite{LDEIM}}
\label{alg:LDEIM}
\end{algorithm}

\section{\label{AP_Computational Cost} Computational Cost}In this Appendix, we perform a computational cost analysis for linear and quadratically nonlinear SPDEs. All of the  computational cost scalings presented in here  exist identically in the DO and BO formulations. To this end,  
let us split the right hand side of the SPDE to a linear and a nonlinear terms: $\mathcal{F}(\bm{v}) =\bm{L}\bm{v} + \bm{N}(\bm{v}) $, where $\bm{L} \in \mathbb{R}^{n \times n}$ and $\bm{N} : \mathbb{R}^n \rightarrow \mathbb{R}^n$ is a nonlinear map. For example, if the SPDE is a one-dimensional Burgers' equation with random initial conditions, $\bm{L}$ is the discrete representation of $\nu \partial^2 ()/\partial x^2$ where $\nu$ is the diffusion coefficient and $\bm{N}(\bm{v})$ is the discrete representation of $-v\partial {v}/\partial {x}$.
\subsection{Homogeneous Linear SPDEs}

First, let us consider a homogeneous linear SPDE where $\bm{N}(\bm{v}) =\bm{0}$, in which case the Burgers' equation reduces to the diffusion equation. For linear equations, the matrix $\bm{F}=\bm{L}\bm{U}\bm{\Sigma}\bm{Y}^T$ does not have to be  computed nor stored explicitly and the DBO equations can be solved in the compressed form. To realize this, consider the right hand side of \cref{eq:DBO_evol_S} and use $\bm{F}=\bm{L}\bm{U}\bm{\Sigma}\bm{Y}^T$. This results in:
\begin{equation*}
    \bm{U}^{T} \bm{W}_{x} \bm{F} \bm{W}_{\xi} \bm{Y} =  \bm{U}^{T} \bm{W}_{x} \bm{L}\bm{U}\bm{\Sigma}\bm{Y}^T \bm{W}_{\xi} \bm{Y} = \bm{U}^{T} \bm{W}_{x} \bm{L}\bm{U}\bm{\Sigma}, 
\end{equation*}
where we have used the orthonormality of the stochastic coefficients $\bm{Y}$ given by \cref{eq:const_Y}. The computational cost of computing $\bm{L}_r \bm{\Sigma} $ is $\mathcal{O}(r^2 n)$, where $\bm{L}_r = \bm{U}^{T} \bm{W}_{x} \bm{L}\bm{U} \in \mathbb{R}^{r \times r}$ is the \emph{reduced linear matrix}. To realize this, note that $\bm{L} \in \mathbb{R}^{n \times n}$ does not have to be stored by utilizing the fact that represents the discretization of spatial derivatives and in fact $\bm{L}$ can be highly sparse. In the case of Burgers' equation: 
\begin{equation}
    \inner{u_i}{\nu \frac{\partial^2 u_j}{\partial x^2}} \approx \big(\bm{L}_r \big)_{ij}, \quad \quad i,j=1, \dots, r,
    \end{equation}
 where $(\bm{L}_r)_{ij}$ is the $(i,j)$ element of  $\bm{L}_r$. Therefore, to compute $\bm{L}_r$ one needs to first compute $\partial^2 u_j/\partial x^2$ for $r$ modes. As an example, if finite difference discretization is used, computing  $\partial^2 u_j/\partial x^2$ is $\mathcal{O}(kn)$, where $k$ is the stencil width of the finite difference scheme. Then the inner product of $u_i$ and $\partial^2 u_j/\partial x^2$ needs to be computed, which is $\mathcal{O}(n)$ and this operation needs to be done  $r^2$ times for $i,j=1,\dots,r$.      The computational cost of  the matrix multiplication $\bm{L}_r \bm{\Sigma} $ is $\mathcal{O}({r^3})$. However, since $r << n$, this cost is negligible in comparison to $\mathcal{O}(n)$.

Let us consider the right hand side of \cref{eq:DBO_evol_U} for a linear SPDE:
\begin{align*}
    \left(\mathbf{I}-\mathbf{U} \mathbf{U}^{T} \mathbf{W}_{x}\right) \mathbf{F} \mathbf{W}_{\xi} \mathbf{Y} \mathbf{\Sigma}^{-1} &= \mathbf{F} \mathbf{W}_{\xi} \mathbf{Y} \mathbf{\Sigma}^{-1} - \mathbf{U} \mathbf{U}^{T} \mathbf{W}_{x} \mathbf{F} \mathbf{W}_{\xi} \mathbf{Y} \mathbf{\Sigma}^{-1}   \\
    &=\bm{L}\bm{U}\bm{\Sigma}\bm{Y}^T \mathbf{W}_{\xi} \mathbf{Y} \mathbf{\Sigma}^{-1} - \mathbf{U} \mathbf{U}^{T} \mathbf{W}_{x} \bm{L}\bm{U}\bm{\Sigma}\bm{Y}^T \mathbf{W}_{\xi} \mathbf{Y} \mathbf{\Sigma}^{-1} \\
     &=\bm{L}\bm{U} - \mathbf{U} \bm{L}_r
\end{align*}
From the analysis of the right hand side of \cref{eq:DBO_evol_S}, it is straightforward to realize that the computational cost of computing $\bm{L}\bm{U} - \mathbf{U} \bm{L}_r$ is also $\mathcal{O}(r^2n)$. Note that  $\bm{L}_r$ needs to be  computed once and it can be utilized in the right hand side of \cref{eq:DBO_evol_S} and \cref{eq:DBO_evol_U}. 

For a homogeneous  linear SPDE it is easy to show that the right hand side of \cref{eq:DBO_evol_Y} is zero:
\begin{align*}
    \left(\mathbf{I}-\mathbf{Y} \mathbf{Y}^{T} \mathbf{W}_{\xi}\right) \mathbf{F}^{T} \mathbf{W}_{x} \mathbf{U} \mathbf{\Sigma}^{-T} &= \left(\mathbf{I}-\mathbf{Y} \mathbf{Y}^{T} \mathbf{W}_{\xi}\right) (\bm{L}\bm{U}\bm{\Sigma}\bm{Y}^T)^{T} \mathbf{W}_{x} \mathbf{U} \mathbf{\Sigma}^{-T}\\
    &=\left(\mathbf{I}-\mathbf{Y} \mathbf{Y}^{T} \mathbf{W}_{\xi}\right) (\bm{Y}\bm{\Sigma}^T\bm{U}^T\bm{L}^T) \mathbf{W}_{x} \mathbf{U} \mathbf{\Sigma}^{-T} \\
    & = \left(\mathbf{Y}-\mathbf{Y} \mathbf{Y}^{T} \mathbf{W}_{\xi} \bm{Y}\right) (\bm{\Sigma}^T\bm{U}^T\bm{L}^T) \mathbf{W}_{x} \mathbf{U} \mathbf{\Sigma}^{-T}\\
    & = \left(\mathbf{Y}-\mathbf{Y} \right) (\bm{\Sigma}^T\bm{U}^T\bm{L}^T) \mathbf{W}_{x} \mathbf{U} \mathbf{\Sigma}^{-T} = \bm{0}
\end{align*}

\begin{remark}
For linear deterministic PDE with random initial conditions, the DBO evolution equations reduce to optimally time-dependent decomposition (OTD) \cite{babaee_minimization_2016, babaee_reduced-order_2017}.
\end{remark}

\subsection{Non-Homogeneous Linear SPDEs}
Now consider the non-homogeneous linear SPDE in the form of $\partial v/\partial t=    \mathcal{L}(v) + g$, where $g(x,t;\omega)$ is a random excitation and $\mathcal{L}$ is a linear differential operator. This equation in the semi-discrete form becomes: $\dot{\bm{V}} = \bm{L}\bm{V} + \bm{G}$, where $\bm{G} \in \mathbb{R}^{n \times s}$ whose columns are random samples of the forcing, the computational complexity of the DBO evolution equation has an additional $\mathcal{O}(rns)$ operation due to the forcing. This can be  observed by investigating the right hand side of \cref{eq:DBO_evol_S}:
\begin{equation*}
    \bm{U}^{T} \bm{W}_{x} \bm{F} \bm{W}_{\xi} \bm{Y} =  \bm{U}^{T} \bm{W}_{x} (\bm{L}\bm{U}\bm{\Sigma}\bm{Y}^T+\bm{G}) \bm{W}_{\xi} \bm{Y} = \bm{L}_r\bm{\Sigma} + \bm{U}^{T} \bm{W}_{x} \bm{G}\bm{W}_{\xi} \bm{Y}, 
\end{equation*}
where 
\begin{equation*}
    \inner{u_i}{\mathbb{E}[gy_j]} \approx \big(\bm{U}^{T} \bm{W}_{x} \bm{G}\bm{W}_{\xi} \bm{Y} \big)_{ij}
\end{equation*}
One can compute $ \bm{U}^{T} \bm{W}_{x} \bm{G}\bm{W}_{\xi} \bm{Y}$ by  first computing  $\bm{G}_{\xi} =\bm{G}\bm{W}_{\xi} \bm{Y}$, which is of order $\mathcal{O}(rns)$ and then computing $ \bm{U}^{T} \bm{W}_{x} \bm{G}_{\xi}$, which is of order $\mathcal{O}(r^2n)$. This term can also be computed by first computing $\bm{G}_{x} = \bm{U}^{T} \bm{W}_{x} \bm{G}$, which is of order $\mathcal{O}(rns)$ and then computing $\bm{G}_{x} \bm{W}_{\xi} \bm{Y}$, which is of order $\mathcal{O}(r^2s)$. Again since $r << s$ and $r <<n$, the overall cost is dominated by $\mathcal{O}(rns)$. The DBO formulation also has $\mathcal{O}(ns)$ memory requirement in cases where $\bm{G}$ is stored in the memory, for example if $\bm{G}$ is time invariant. If memory limitations does not allow that, then rows or columns of $\bm{G}$   must be computed one (or few) at a time.  

\subsection{Nonlinear SPDEs}
The computational cost of computing the right hand side terms of \cref{eq:DBO_evol_S,eq:DBO_evol_U,eq:DBO_evol_Y} for nonlinear terms scale with $\mathcal{O}(sn)$ since the nonlinear term $\bm{N}(\bm{V})$ must be computed for all columns of $\bm{V}$, i.e., $\bm{N}(\bm{V})=[\bm{N}(\bm{v}_1), \bm{N}(\bm{v}_2), \dots , \bm{N}(\bm{v}_s)]$. For quadratic nonlinearities, e.g., the Burgers' equation,  it is possible  to compute the projection of  $\bm{N}(\bm{V})$ onto spatial and stochastic bases in a compressed form, i.e., by  not forming the matrix $\bm{N}(\bm{V})$.  To see this, first we note that: 
\begin{equation*}
\inner{u_m}{-\mathbb{E}[(u_i\Sigma_{ij}y_j)(\frac{\partial u_{i'}}{\partial x}\Sigma_{i'j'}y_{j'}) y_n]} \approx \big(\bm{U}^T\bm{W}_x \bm{N}(\bm{U}\bm{\Sigma}\bm{Y}^T)\bm{W}_{\xi}\bm{Y}  \big)_{mn} 
\end{equation*}
where  
 \begin{equation*}
\bm{N}(\bm{U}\bm{\Sigma}\bm{Y}^T) = -\big(\bm{u}_i\odot (\bm{D}_x \bm{u}_{i'}) \big) \bm{\Sigma}_{ij}\bm{\Sigma}_{i'j'}\big(\bm{y}^T_j\odot  \bm{y}^T_{j'}\big). 
\end{equation*}
Here, $\odot$ represents an element-wise product between two vectors and $\bm{D}_x \in \mathbb{R}^{n \times n}$ is the discrete representation of $\partial ()/\partial x$ and the repeated indices imply summation over those indices. The computational cost of computing  $\bm{u}_i\odot (\bm{D}_x \bm{u}_{i'})$ for $i,i'=1,\dots, r$ is $\mathcal{O}(r^2 n)$,  since the computational cost of computing  $\bm{u}_i\odot (\bm{D}_x \bm{u}_{i'})$ for each pair of $(i,i')$ is $\mathcal{O}(n)$. Similarly, the computational cost of computing $\bm{y}^T_j\odot  \bm{y}^T_{j'}$ is $\mathcal{O}(r^2s)$. Therefore, for quadratic nonlinearity, it is possible to reduce  $\mathcal{O}(sn)$ to $\mathcal{O}(r^2(n+s))$. However, this can be achieved in an \emph{intrusive} manner, i.e., by replacing the DBO expansion into the nonlinear form and derive and implement  the resulting nonlinear terms. It is straightforward to show that for polynomial nonlinearity of order $m$, the same approach results in the computational complexity of $\mathcal{O}(r^m n)+\mathcal{O}(r^ms)$. Therefore, the cost increases exponentially fast with $m$.  Also, as $m$ increases $(m>2)$, deriving and implementing the nonlinear expansion of the DBO decomposition can become  overwhelming due to the highly intrusive nature of this approach, which generates exponentially larger number of  terms as $m$ increases.  For non-polynomial nonlinearity, e.g., exponential and rational nonlinearities,  it is not possible to avoid  $\mathcal{O}(sn)$ cost because the nonlinear expansion,  for example $\exp(\bm{U}\bm{\Sigma}\bm{Y}^T)$, requires infinitely many terms.

\bibliographystyle{elsarticle-num}
\biboptions{sort&compress}
\bibliography{ms}

\end{document}